\documentclass[10pt]{article}

\usepackage{amsmath, amssymb, amsfonts, amsthm, graphicx}

\DeclareMathOperator{\Area}{Area}
\DeclareMathOperator{\bd}{bd}
\DeclareMathOperator{\conv}{conv}
\DeclareMathOperator{\pos}{pos}
\DeclareMathOperator{\interior}{int}
\DeclareMathOperator{\rel}{rel}
\newcommand{\normal}{\dashv}
\newcommand{\notnormal}{\not\,\dashv}
\newcommand{\Line}[1]{\langle#1\rangle}
\newcommand{\Ray}[1]{[#1\rangle}
\newcommand{\Segment}[1]{[#1]}
\newcommand{\epsi}{\varepsilon}
\newcommand{\norm}[1]{\|#1\|}
\newcommand{\un}[1]{\widehat{\!#1}}

\newcommand{\R}{\mathbf{R}}
\newcommand{\Z}{\mathbf{Z}}
\newcommand{\N}{\mathbf{N}}
\newcommand{\vect}[1]{\mbox{\itshape\bfseries#1}}
\newcommand{\va}{\vect{a}}
\newcommand{\vb}{\vect{b}}
\newcommand{\vc}{\vect{c}}
\newcommand{\vd}{\vect{d}}
\newcommand{\ve}{\vect{e}}
\newcommand{\vf}{\vect{f}}
\newcommand{\vx}{\vect{x}}
\newcommand{\vy}{\vect{y}}
\newcommand{\vz}{\vect{z}}

\newcommand{\vo}{\vect{o}}
\newcommand{\vp}{\vect{p}}
\newcommand{\vq}{\vect{q}}
\newcommand{\vr}{\vect{r}}
\newcommand{\vu}{\vect{u}}
\newcommand{\vv}{\vect{v}}
\newcommand{\vw}{\vect{w}}

\newcommand{\curve}[1]{{\mathcal#1}}

\newcommand{\abs}[1]{|#1|}

\newcommand{\myangle}{\sphericalangle}
\newcommand{\length}[1]{|#1|}

\theoremstyle{plain}
\newtheorem{theorem}{Theorem}
\newtheorem{lemma}[theorem]{Lemma}
\newtheorem{corollary}[theorem]{Corollary}
\newtheorem{proposition}[theorem]{Proposition}

\theoremstyle{definition}

\begin{document}

\bibliographystyle{amsplain}

\title{The geometry of Minkowski spaces --- a survey. Part I}
\author{Horst Martini \\
        Fakult\"at f\"ur Mathematik \\
        Technische Universit\"at Chemnitz \\
        D-01097 Chemnitz \\ Germany \\
        E-mail: \texttt{martini@mathematik.tu-chemnitz.de} \and 
        Konrad J. Swanepoel\thanks{Research supported by a grant from a cooperation between the Deutsche For\-schungs\-gemein\-schaft in Germany and the National Research Foundation in South Africa}\\ 
        Department of Mathematics and Applied Mathematics \\ 
        University of Pretoria, Pretoria 0002 \\ South Africa \\
        E-mail: \texttt{konrad@math.up.ac.za} \and 
        Gunter Wei\ss \\
        Institut f\"ur Geometrie \\
        Technische Universit\"at Dresden \\
        D-01062 Dresden \\ Germany \\
        E-mail: \texttt{weiss@math.tu-dresden.de}}
\date{}
\maketitle

\begin{abstract}
We survey elementary results in Minkowski spaces (i.e.\ finite dimensional Banach spaces) that deserve to be collected together, and give simple proofs for some of them.
We place special emphasis on planar results.
Many of these results have often been rediscovered as lemmas to other results.
In Part I we cover the following topics:
The triangle inequality and consequences such as the monotonicity lemma, geometric characterizations of strict convexity, normality (Birkhoff orthogonality), conjugate diameters and Radon curves, equilateral triangles and the affine regular hexagon construction, equilateral sets, circles: intersection, circumscribed, characterizations, circumference and area, inscribed equilateral polygons.
\end{abstract}

\tableofcontents

\newpage
\section{Introduction}
This paper is a survey of basic results on the geometry of finite dimensional normed linear spaces, which we, following Thompson \cite{MR97f:52001}, call \emph{Minkowski geometry}.
We place special emphasis on planar results, since there are many planar results which are simple and elementary, or at least considered to be so, which are often needed in work on Minkowski geometry, often rediscovered, and not often proved.
They come from various parts of mathematics: discrete geometry and geometry of numbers (the historical origin of Minkowski geometry), convex geometry, functional analysis, and lately also from optimization, theoretical computer science and combinatorics.
Our motivation for this survey is that the book of Thompson \cite{MR97f:52001} does not cover many of these results; neither does the Handbook of Convex Geometry \cite{MR94e:52001} (as admitted in its introduction).
Previous surveys of this type, such as Petty \cite{MR18:760e} and Yaglom \cite{MR58:30739}, are old and not easily accessible.
In many of these results, the proofs, once written out, have a tendency to be messy.
In this regard Sch\"affer \cite{MR36:1959} says (in 1967)
\begin{quote}
The amazing amount of underbrush that has to be cleared away, \dots, indicates, to this author at least, that the geometry of finite-dimensional convex sets is still quite imperfectly known.
\end{quote}
We discuss many of these elementary results, and give proofs whenever the original proofs are in journals that are difficult to find, or if we have simpler proofs than those in the literature.
In each case we then survey the extensions to higher dimensions.
There are many topics that are not included here, but which will appear in a continuation of this survey, such as $d$-segments and $d$-convexity, characterizations of smoothness, various notions of angle measures, diametrically maximal sets, sets of constant width, packing and covering of unit balls, various discrete inequalities, Hadwiger numbers, bisectors, Erd\H os-type problems, approximation theory (Chebyshev sets), isometries (the Mazur-Ulam theorem and its relatives, Beckman-Quarles type theorems, Banach-Mazur distance), parameters of functional analysis, applications to discrete optimization problems such as minimum spanning trees and Steiner minimum trees, the Fermat-Torricelli problem, etc.
We do not consider the local theory of Banach spaces itself \cite{MR90k:46039, MR91d:52005} nor differential and integral geometry of Minkowski spaces, for which the main reference is \cite{MR97f:52001}.
We have here mainly concentrated on the triangle inequality and its consequences such as the monotonicity lemma; geometric characterizations of strict convexity and smoothness; normality (Birkhoff orthogonality), conjugate diameters and Radon curves; equilateral triangles and the affine regular hexagon construction; equilateral sets; various aspects of circles and spheres: intersection of circles, circumscribed spheres, characterizations of spheres, circumference and area of the unit circle, and equilateral polygons inscribed in the unit circle.

Also, we do not consider infinite dimensional spaces, and for the sake of simplicity, we only consider Minkowski spaces with a symmetric norm ($\norm{\vx} = \norm{-\vx}$) satisfying the triangle inequality.

\section{The subject and its origin}\label{definitions}
As is well known, the axioms of Minkowski spaces were introduced by Minkowski \cite{Minkowski0}, in connection with problems from number theory.
However, it seems that the earliest reference to non-Euclidean geometry in the sense of Minkowski Geometry was made by Riemann in his Habilitationsvortrag \cite{Riemann}, where he mentioned the $\ell_4$-norm.
See \'Alvarez \cite{Alvarez} for a discussion of Riemann's remarks on non-Euclidean norms.
Hilbert \cite{Hilbert} in his famous lecture delivered before the International Congress of Mathematicians in 1900 gives a description of Minkowski Geometry in his fourth problem.

Two important and neglected early papers, considering Minkowski Geometry from a geometric (as opposed to analytic) point of view are \cite{Golab} and \cite{GolabHarlen}.
Minkowski Geometry was studied, especially by Busemann \cite{MR17:779a}, in order to throw more light on Finsler Geometry, introduced by Finsler \cite{Finsler}.
See also \cite{MR21:4462}.
For recent developments in Finsler Geometry, see \'Alvarez \cite{Alvarez}.
Closely related is the subject of distance geometry, cf.\ the work of Menger and his school, summarized by Menger \cite{Menger} and Blumenthal \cite{MR42:3678} (see also \cite{MR42:8370}), as well as the works of Alexandrov \cite{MR10:619c} and Rinow \cite{MR23:A1290}.

In Functional Analysis, although concentrating from its outset almost exclusively on infinite dimensional spaces \cite{Banach, MR58:2112}, various fine geometric properties of finite dimensional spaces play an important role in the so-called local theory of Banach spaces \cite{MR87m:46038, MR90k:46039, MR91d:52005, MR95b:46012}.
For infinite dimensional Banach space geometry, see \cite{MR49:9588, MR57:1079, MR58:17766, MR88f:46021}.
Also, characterizations of inner product spaces lead to interesting geometry of finite dimensional spaces \cite{MR88m:46001}.

Recently, in Operations Research and VLSI design, various norms have started playing an important role, especially $\ell_p$-norms and polygonal norms; see e.g.\ \cite{MR1358610, MR99i:05062}.

For physical interpretations of Minkowski spaces, see e.g.\ \cite{Laugwitz, MR93h:53012, MR95i:58051}.
The two most common equivalent definitions of a Minkowski space are by giving axioms for a norm $\norm{\cdot}:V\to\R$, namely
\[ \begin{array}{rcll}
\norm{\vx} & \geq & 0, \\
\norm{\vx} & = & 0 \mbox{ iff }  \vx=\vo & \mbox{ (positive definiteness)},\\
\norm{\lambda\vx} & = & \abs{\lambda}\,\norm{\vx} & \mbox{ (symmetry)},\\
\norm{\vx+\vy}& \leq & \norm{\vx}+\norm{\vy} & \mbox{ (triangle inequality)},
\end{array} \]
where $V$ is the underlying finite dimensional real vector space,
or by giving axioms for the unit ball $B\subset V$, namely
\[ \begin{array}{l}
B \mbox{ is bounded,} \\
B \mbox{ has a non-empty interior,} \\
B \mbox{ is centrally symmetric,} \\
B \mbox{ is convex.}
\end{array} \]

In the sequel, $M=M_d$ will be an arbitrary $d$-dimensional Minkowski space with norm $\norm{\cdot}$, unit ball $B=B(M)$, and unit sphere $S=S(M)$.
A two-dimensional Minkowski space will be referred to as a \emph{Minkowski plane}, its unit ball called the \emph{unit disc} and its unit sphere the \emph{unit circle}.

We identify isometric spaces; by the Mazur-Ulam theorem \cite{MazurUlam} two Minkowski spaces are isometric iff their unit balls are affinely equivalent.
The Minkowski plane with a parallelogram as unit disc will be called the \emph{rectilinear} or \emph{taxicab plane}.
We denote vectors by $\vx,\vy,\dots$, the normalization of $\vx\neq \vo$ by $\un{\vx}:= \frac{1}{\norm{\vx}}\vx$, the closed (straight line) segment from $\va$ to $\vb$ by $\Segment{\va\vb}$, the triangle with vertices $\vp,\vq,\vr$ by $\triangle\vp\vq\vr$, the line through $\va$ and $\vb$ by $\Line{\va\vb}$, the ray with origin $\va$ passing through $\vb$ by $\Ray{\va\vb}$, curves parametrized by $[a,b]$ by $\gamma$, the length of a segment $\Segment{\va\vb}$ (in the norm) by $\length{\va\vb} := \norm{\vb-\va}$, and the length of a curve by $\length{\gamma}$, which we define in the elementary way (i.e.\ without integrals)
\[\length{\gamma}:=\sup\{\sum_{i=1}^n\length{\va_i\va_{i+1}} : n\in\N, \va_i=\gamma(t_i), a=t_0<t_1<\dots<t_{n-1}=b\}.\]
The unit circle $S$ of a Minkowski plane $M$, parametrized as a curve, has a length $\ell(S)$, called its \emph{circumference}.
We use the following notation: $\Pi(M) := \ell(S)/2$.
A \emph{metric segment} is a curve which is isometric to a closed segment of the real line,
a \emph{metric line} is a curve parametrized by $\R$ and isometric to the real line, and a \emph{geodesic} (a notion due to Busemann \cite{MR4:109e}; see also \cite[p.~32]{MR17:779a}) is a curve that is locally a metric segment, i.e., each point of the curve has a closed neighbourhood that is a metric segment.

We use the standard facts about duals, denoting functionals by $\phi,\psi$ etc.
An important fact is that all Minkowski spaces are reflexive Banach spaces.
We denote the convex hull of a set $S$ by $\conv S$, the positive hull by $\pos S$, the interior by $\interior S$, the relative interior by $\rel\interior S$, and the boundary by $\bd S$.
A (positive) \emph{homothet} of a set $A$ is a set of the form $\lambda A + \vv$, where $\lambda>0$ and $\vv$ are arbitrary.
Thus we allow translates as special homothets (with the point of homothety at infinity).

We say that $\vx\neq \vo$ is normal to $\vy\neq\vo$, denoted $\vx\normal\vy$, if $\norm{\vx}\leq\norm{\vx+\lambda\vy}$ for all $\lambda\in\R$, or equivalently, if the unit ball is supported at $\un{\vx}$ by a line parallel to $\vy$.
See Section~\ref{normalsection} for a discussion of this notion.

The difference body of a convex body $C$ is $\frac{1}{2}(C-C)$.
Given any convex body $C$, its difference body is centrally symmetric, hence induces a norm for which $C$ has constant width.
A \emph{face} of a convex body $C$ is a subset $F$ of $C$ such that, whenever the relative interior of some segment $\Segment{\va\vb}\subseteq C$ intersects $S$, then $\Segment{\va\vb}\subseteq S$.
It is a well-known fact from convex geometry that any point of a convex body is contained in the relative interior of a unique face (see \cite[Theorem 2.1.2]{MR94d:52007}).
A set of points $S$ in a finite dimensional vector space is \emph{antipodal} if for any distinct $\vx,\vy\in S$ there exist parallel (distinct) hyperplanes supporting $\conv S$ at $\vx$ and $\vy$.
The set $S$ is \emph{strictly antipodal} if it is furthermore possible to find hyperplanes containing no other point of $\conv S$.

\section{The triangle inequality}\label{triineqsect}
We only consider the triangle inequality and some of its simple consequences.
We omit altogether a discussion of hypermetric inequalities, piecewise
linear inequalities, duals of zonoids, etc., since this is
well-treated in e.g.\ the Handbook of Convex Geometry \cite{MR94e:52001}.
We also postpone discussion of Fermat-Torricelli points, sums of
distances, Chebyshev centres, vector balancing and discrepancy theory, etc., to part II.
The books of Bottema et al.\ \cite{MR41:7537} and Mitrinovi{\'c} et al.\ \cite{MR91k:52014} have many inequalities for triangles, which in fact are true for all metric spaces, since they are purely algebraic consequences of the triangle inequality (e.g.\ inequalities 1.1--1.24 in \cite{MR41:7537}).

We now discuss the triangle inequality and metric lines.

When one defines a distance function using only a star-shaped body, then the triangle inequality is equivalent to the convexity of unit ball, and the strict triangle inequality (equality only for collinear points) is equivalent to strict convexity.
This was already proved by Minkowski \cite[Kapitel 1, Paragraph 18]{Minkowski1}.

In the following two propositions we analyse the triangle inequality further.
The first, characterizing specific instances of equality in the triangle inequality, has been noticed by Alt \cite{Alt}, and also by
Go{\l}ab and H\"arlen \cite{GolabHarlen}.

\begin{proposition}\label{trieq}
For any distinct $\va,\vb,\vc$ in a Minkowski space, $\length{\va\vc}\leq\length{\va\vb}+\length{\vb\vc}$, with equality iff $\Segment{\vx\vy}\subseteq S$, where $\vx=\un{\vb-\va}$ and $\vx=\un{\vc-\vb}$.
\end{proposition}

\proof
Let $\alpha=\length{\va\vb}, \beta=\length{\vb\vc},\gamma=\length{\va\vc}$.
Letting $\vz=\un{\vc-\va}$, we have $\vz=(\alpha/\gamma)\vx+(\beta/\gamma)\vy$.

Suppose now $\gamma=\alpha+\beta$.
Then we see that the above representation of $\vz$ is a convex combination of $\vx$ and $\vy$, with $\vz$ strictly between $\vx$ and $\vy$.
Thus we have three collinear vectors on the boundary of the unit ball, and it follows that their convex hull is also on the boundary, i.e.,
$\Segment{\vx\vy}$ is a segment on the unit sphere.

Conversely, suppose that $\Segment{\vx\vy}$ is a segment on the unit sphere.
Let $\mu=\gamma/(\alpha+\beta)$.
Thus we have the convex combination $\mu\vz= \mu(\alpha/\gamma)\vx+\mu(\beta/\gamma)\vy$, i.e., $\mu\vz$ is on the segment $\Segment{\vx\vy}$.
It follows that $\mu\vz$ is a unit vector, $\mu=1$, and $\gamma=\alpha+\beta$.
\qed

The above proposition has also been observed recently by Nitka and Wiatrowska \cite{MR39:2070}, who furthermore make the remark that, given any three points in $d$-dimensional vector space, there is a norm such that the three points give equality in the triangle inequality.
One merely has to give a unit ball with the right segment on its boundary.

Go{\l}ab and H\"arlen \cite{GolabHarlen} made a very complete analysis of the triangle inequality.
The following proposition, due to them, states that the extreme points of the unit ball coincide with the directions of strict inequality in the triangle inequality, and with the directions of unique \emph{metric segments}, i.e., curves $\gamma$ from $\va$ to $\vb$ such that $\ell(\gamma)=\length{\va\vb}$.
This result has also been observed recently by Toranzos \cite{MR46:4171}.

\begin{proposition}
Let $\vv$ be a unit vector.
Then the following are equivalent.
\begin{enumerate}
\item $\vv$ is an extreme point of the unit ball.
\item For any distinct $\va,\vb,\vc$ such that $\vv=\un{\vc-\va}$ we have $\length{\va\vc}=\length{\va\vb}+\length{\vb\vc}$ iff $\va,\vb,\vc$ are collinear.
\item The segment $\Segment{\va\vb}$ is the unique metric segment joining $\va$ to $\vb$.
\end{enumerate}
\end{proposition}

Metric segments and lines in a more general context have also been investigated by Menger \cite{Menger}.

It is also now simple to characterize metric segments in general.
Note that a curve from $\va$ to $\vb$ is a metric segment iff it has length $\length{\va\vb}$.
See Section~\ref{definitions} for the definition of a face of a convex body.

\begin{proposition}\label{msegchar}
A curve $\gamma$ from $\va$ to $\vb$ is a metric segment iff each directed chord of $\gamma$ is in a direction contained in the unique face of the unit ball containing $\un{\vb-\va}$ in its relative interior.
\end{proposition}

\begin{figure}
\begin{center}
\includegraphics{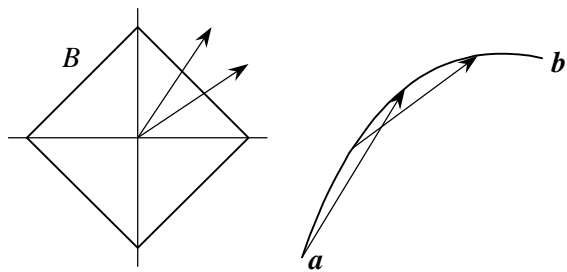}
\end{center}
\caption{Characterizing metric segments by their chords}
\end{figure}

There are at least $2d$ directions in which there are unique metric segments (Go{\l}ab and H\"arlen show at least $d+1$ in the case of a non-symmetric norm), and more if the unit ball is not a cross-polytope.
This can be seen from the higher-dimensional generalization of Proposition~\ref{conjdiam} (in Section~\ref{conjdiamsection} below).

\begin{figure}
\begin{center}
\includegraphics{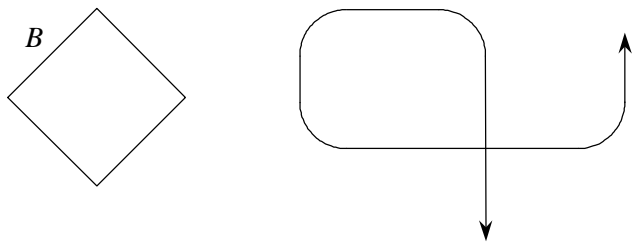}
\end{center}
\caption{A geodesic that is not a metric segment}\label{fig2}
\end{figure}

It was noticed by Szenthe \cite{MR23:A2779} that in certain Minkowski spaces there exist geodesics which are not metric segments.
For example, using Proposition~\ref{msegchar}, the curve in Figure~\ref{fig2} is easily seen to be a geodesic.

Szenthe gave the following characterization of Minkowski spaces which have geodesics that are not metric segments.

\begin{theorem}
A Minkowski space has geodesics which are not metric segments iff the unit sphere contains two segments with a common endpoint that are not contained in the same supporting hyperplane.
\end{theorem}

We now consider special consequences of the triangle inequality in Minkowski planes.

\subsection{The triangle inequality in Minkowski planes}

The following two lemmas also hold for higher dimensional spaces, but the argument is two-dimensional.
They are mentioned in e.g.\ \cite{MR93e:52038}.

\begin{lemma}\label{maxineq}
If $\vw$ is strictly between $\vy$ and $\vz$, then $\length{\vx\vw}\leq\max\{\length{\vx\vy},\length{\vx\vz}\}$, with equality iff $\length{\vx\vw}=\length{\vx\vy}=\length{\vx\vz}$.
In the case of equality, $\length{\vx\vw}$ is the shortest distance from $\vx$ to the line $\Line{\vw\vz}$, and $\length{\vx\vw}=\length{\vx\vv}$ for all $\vv\in\Segment{\vy\vz}$.
\end{lemma}

\begin{figure}
\begin{minipage}{6.5cm}
\begin{center}
\bigskip
\includegraphics{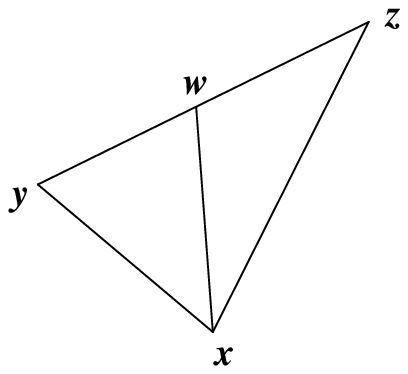}
\end{center}
\caption{$\length{\vx\vw}\leq\max\{\length{\vx\vy},\length{\vx\vz}\}$}
\end{minipage}
\hfill
\begin{minipage}{5.5cm}
\begin{center}
\includegraphics{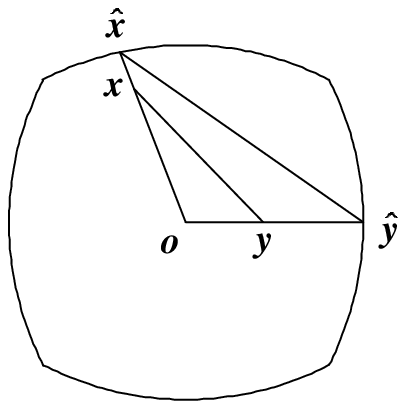}
\end{center}
\caption{If $\length{\vx\vy}\geq 1$, then $\length{\un{\vx}\un{\vy}}\geq 1$}
\end{minipage}
\end{figure}

\proof
Writing $\vw=\lambda\vy+(1-\lambda)\vz$ for some $0<\lambda< 1$ we have the following sequence of inequalities:
\begin{eqnarray}
\length{\vx\vw} & = & \norm{\vx-(\lambda\vy+(1-\lambda)\vz)} \nonumber \\
& \leq & \lambda\length{\vx\vy} + (1-\lambda)\length{\vx\vz} \label{maxineq1} \\
& \leq & \max(\length{\vx\vy},\length{\vx\vz}). \label{maxineq2}
\end{eqnarray}

Suppose now that we have equality.
Equality in \eqref{maxineq2} gives $\length{\vx\vy}=\length{\vx\vz}$, hence also $=\length{\vx\vw}$.
Thus $\vy,\vw,\vz$ are on the circle with centre $\vx$ and radius $\length{\vx\vw}$.
It follows that $\Segment{\vy\vz}$ is on that circle.
Thus $\length{\vx\vw}=\length{\vx\vv}$ for all $\vv\in\Segment{\vy\vz}$.

It remains to prove that this is the minimum distance.
Let $\vp$ be such that $\vz$ is between $\vy$ and $\vp$.
Suppose $\length{\vx\vp}<\length{\vx\vy}$.
Applying the above inequality to $\vy, \vz$ and $\vp$, we obtain that there is in fact equality $\length{\vx\vz}=\max(\length{\vx\vy},\length{\vx\vp})$, hence $\length{\vx\vy}=\length{\vx\vp}$, a contradiction.
Thus $\length{\vx\vp}\geq\length{\vx\vy}$.

We get the same inequality for $\vp$ such that $\vy$ is between $\vp$ and $\vz$.
\qed

\begin{lemma}\label{packingineq}
Let $\vx,\vy\neq\vo$ be contained in the unit ball.
If $\length{\vx\vy}\geq 1$, then $\length{\un{\vx}\un{\vy}}\geq 1$.
\end{lemma}

\proof
By multiplying $\vx,\vy$ by a suitable scalar $\geq 1$, we may assume without loss of generality that $\norm{\vx}=1\geq\norm{\vy}$.
By Lemma~\ref{maxineq} it follows that $\length{\vx\vy}\leq\max(\length{\vx\vo},\length{\vx\un{\vy}})$.
If $\length{\vx\un{\vy}}<1$, then, since $\length{\vx\vy}\geq 1$, we have equality in the above inequality, and it follows from Lemma~\ref{maxineq} that $\length{\vx\un{\vy}}=1$, a contradiction.
Hence $\length{\vx\un{\vy}}\geq 1$.
\qed

We note that Lemma~\ref{maxineq} also immediately gives that the diameter of a convex body in a Minkowski space is attained at extreme points.

\begin{proposition}\label{quadlemma}
In a convex quadrilateral $\va\vb\vc\vd$ in a Minkowski plane, the sum of the diagonals is at least the sum of two opposite sides, i.e.,
\begin{equation}\label{eq1}
\length{\va\vc}+\length{\vb\vd}\geq\length{\va\vb}+\length{\vc\vd}
\end{equation}
and
\begin{equation*}
\length{\va\vc}+\length{\vb\vd}\geq\length{\vb\vc}+\length{\va\vd},
\end{equation*}
with equality in \eqref{eq1}, say, iff $\Segment{\un{\vc-\va}\;\;\un{\vb-\vd}}\subseteq S$.
\end{proposition}

\begin{figure}
\begin{minipage}[b]{5.5cm}
\begin{center}
\includegraphics{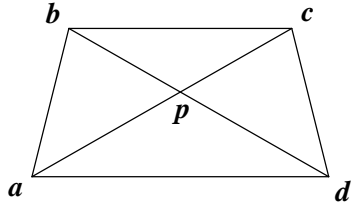}
\end{center}
\caption{The sum of the diagonals is at least the sum of opposite sides}
\end{minipage}
\hfill
\begin{minipage}[b]{6cm}
\begin{center}
\includegraphics{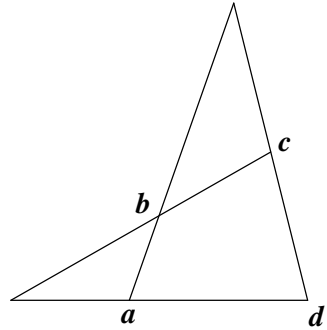}
\end{center}
\caption{Counterpart to Proposition~\ref{quadlemma}}
\end{minipage}
\end{figure}

\begin{corollary}
In a convex quadrilateral $\va\vb\vc\vd$ in a Minkowski plane, twice the sum of the diagonals is at least the sum of the sides, with equality iff the plane is rectilinear.
\end{corollary}

Both the proposition and the corollary above follow from Proposition~\ref{trieq} using the intersection point of the diagonals.

It follows that if the plane is not rectilinear one cannot have the diagonals as the two shortest distances between the four points, and if the plane is strictly convex one cannot have two opposite sides as the two largest distances between the four points (mentioned by Gr\"unbaum and Kelly \cite[Theorem 1]{MR39:6180}).
See Brass \cite{MR97c:52036} for combinatorial consequences.
Also note that the Dowker-type results \cite{MR5:153m} for perimeter of Moln{\'a}r \cite{MR17:1235c} (see \cite[Chapter 2]{MR96j:52001} for the Euclidean case) are immediately extendable to Minkowski planes, using elementary inequalities still valid for Minkowski planes: for inscribed $n$-gons one needs Proposition~\ref{quadlemma}, and for circumscribed $n$-gons one needs 

\begin{proposition}
In a convex quadrilateral $\va\vb\vc\vd$ in a Minkowski plane, if $\Line{\va\vd}$ and $\Line{\vb\vc}$ intersect on the side of $\Line{\va\vb}$ opposite $\vc$ and $\vd$, and $\Line{\va\vb}$ and $\Line{\vc\vd}$ intersect on the side of $\Line{\vb\vc}$ opposite $\va$ and $\vd$, then
\[\length{\va\vb}+\length{\vb\vc}\leq\length{\va\vd}+\length{\vd\vc}.\]
\end{proposition}

The remainder of the proof (see \cite{MR96j:52001}) is then easily adapted for Minkowski planes, and we obtain
\begin{theorem}
Let $C$ be a convex disc in  a Minkowski plane $M$.
For each $n\geq 3$, let $Q_n$ ($q_n$) be an $n$-gon of minimum (maximum) perimeter circumscribed about $C$ (inscribed in $C$).
Then
\[ \length{Q_n} \leq \frac{\length{Q_{n-1}}+\length{Q_{n+1}}}{2}\mbox{ and }
   \length{q_n} \geq \frac{\length{q_{n-1}}+\length{q_{n+1}}}{2} \]
for all $n\geq 4$.
\end{theorem}

\subsection{Strict convexity}\label{scsect}
A Minkowski space is \emph{strictly convex} if $\norm{\vx+\vy}=\norm{\vx}+\norm{\vy}$ implies that $\vx$ and $\vy$ are linearly dependent, i.e., there is equality in the triangle inequality only in the trivial cases.

The unit ball of a Minkowski space is \emph{rotund} if its boundary does not contain any straight line segment, i.e., each boundary point is an extreme point.

There are various strengthenings of the definition of strict convexity, such as uniform convexity, but these stronger concepts are usually only stronger in the infinite dimensional context; see the survey by Cudia \cite{MR27:5106}.
We now discuss a few properties that are equivalent to strict convexity of a Minkowski space.
We only consider geometric characterizations and ignore the many characterizations in terms of operators, duality maps and semi-inner products such as those given in \cite{MR41:5943, MR54:5806, MR84b:46021}.

Note that strict convexity is a two-dimensional notion, i.e., a space is strictly convex iff each of its two-dimensional subspaces is strictly convex.
The following characterizations of strict convexity for Minkowski spaces are essentially folklore results from convex geometry (see Minkowski \cite[Kapitel 1, Paragraph 18]{Minkowski1} and Day \cite[p.~144]{MR49:9588}).
See also Bumcrot for a further discussion \cite{MR38:3714}.

The following are equivalent to strict convexity of a Minkowski space:
\begin{enumerate} 
\item every boundary point is an extreme point (exposed point),
\item metric segments are always straight line segments,
\item the unit ball is rotund,
\item a linear functional has at most one maximum on the unit ball,
\item any supporting hyperplane of the unit ball touches the unit ball in at most one boundary point,
\item supporting hyperplanes at distinct points of the boundary of the unit ball are distinct.
\end{enumerate}

James \cite{MR9:42c} gave a few characterizations in terms of properties of normality:
A Minkowski space is strictly convex iff normality is left unique, i.e., for all $\vx_1,\vx_2\neq \vo$, if $\vx_1\normal\vy$ and $\vx_2\normal\vy$ for some $\vy\neq\vo$, then $\vx_1=\vx_2$.

In terms of nearest points there are the following equivalences:
\begin{enumerate}
\item For each point $\vp$ and each convex set $C$ there is at most one point in $C$ that is nearest to $\vp$.
\item For each point $\vp$ and each closed convex set $C$ there is exactly one point in $C$ that is nearest to $\vp$.
\item For each point $\vp$ and each one-dimensional subspace $L$ there is exactly one point in $L$ that is nearest to $\vp$.
\item For each point $\vp$ and each metric line $\gamma$ there is exactly one point in $\gamma$ that is nearest to $\vp$ (Andalafte and Valentine \cite{MR47:2499}).
\end{enumerate}
Singer \cite{MR24:A1629} gives generalizations of the above to properties weaker than strict convexity.

Andalafte and Valentine \cite{MR47:2499} note the following two characterizations.

\begin{proposition}
The following properties are equivalent to strict convexity of a Minkowski space.
\begin{enumerate}
\item A line and a sphere intersect in at most two points.
\item The distance from a fixed point to a variable point on a line is strictly unimodal.
\end{enumerate}
\end{proposition}

\proof
1 immediately follows from 2.
In 2 unimodality follows immediately from convexity of the mentioned function.
Strict unimodality is then equivalent to a unique minimum, which is number 3 of the above-mentioned characterizations in terms of nearest points.
\qed

Note that the above-mentioned unimodality is formulated in an elementary way in Lemma~\ref{maxineq}.

Valentine \cite{MR84m:46023} gives the following characterizations.

\begin{proposition}
The following properties are equivalent to strict convexity of a Minkowski space.
\begin{enumerate}
\item The diagonals of a metric parallelogram \textup{(}i.e., a planar quadrilateral with opposite sides of equal length\textup{)} bisect each other.
\item There is no \emph{ramification point}, i.e., a point $\vx$ such that there are distinct points $\va,\vb,\vc$ such that $\length{\va\vx}+\length{\vx\vb}=\length{\va\vb}$ and $\length{\va\vx}+\length{\vx\vc}=\length{\va\vc}$.
\end{enumerate}
\end{proposition}

There are statements in Euclidean geometry about lengths that characterize equality in the triangle inequality.
These then also lead to characterizations of strict convexity.
As an example, the Heron formula for the area of a triangle in the Euclidean plane gives $0$ exactly when the sum of two sides equals the third side.
This then gives a characterization that has essentially the same geometric content as the definition (see Reda \cite{MR52:4244} and Diminnie and White \cite{MR57:13452}).

See Diminnie and White \cite{MR82k:52004} for characterizations of strict convexity in terms of various postulates for the metric betweenness relation.

We now turn to characterizations in Minkowski planes.
Holub \cite{MR52:1263} gives two interesting characterizations of strict convexity of Minkowski planes.
We begin by stating a lemma that was needed by Benz for a Beckman-Quarles type theorem \cite{MR88b:51022, MR88j:46020}, and was also discussed in \cite{MR90a:46025}.

\begin{lemma}
In a Minkowski plane,
\begin{enumerate}
\item any vector of norm $<1$ is the midpoint of some chord of the unit circle,
\item any vector of norm $<2$ is the sum of two unit vectors,
\item any two unit circles with centres at distance $<2$ intersect.
\end{enumerate}
\end{lemma}

\proof
It is easily seen that the three statements are equivalent, while the last is intuitively obvious, and follows from a simple continuity argument.
\qed

Holub's one characterization is that a Minkowski plane is strictly convex iff any vector of norm $<1$ is the midpoint of \emph{at most} one chord of the unit circle.
In fact, it is easily seen that we may replace ``$<1$'' by ``$<\epsi$ for some $\epsi>0$''.
Also, it must have been noticed very early that two circles in a strictly convex plane intersect in at most two points.
Two early references are Mayer \cite{Mayer}, who also assumes smoothness, and Buter \cite{Buter}, who drops the assumption of smoothness.
Mayer also remarks that this is equivalent to the fact that any three points are contained in at most one unit circle.
Valentine \cite{MR84m:46023} uses a different geometric approach to show that strict convexity is characterized by the fact that two unit circles intersect in at most one point on each side of the line passing through the two centres.
His approach is related to the Monotonicity Lemma (see Section~\ref{monotonsubsect}).

By these remarks we now have the following list of equivalences.

\begin{proposition}
Each of the following statements are equivalent to strict convexity of a Minkowski plane.
\begin{enumerate}
\item Any vector of norm $<1$ is the midpoint of a unique \textup{(}at most one\textup{)} chord of the unit circle,
\item Any vector of norm $<2$ is the sum of two unit vectors in a unique \textup{(}at most one\textup{)} way,
\item Any two unit circles with centres at distance $<2$ intersect in exactly \textup{(}at most\textup{)} two points,
\item There is an $\epsi>0$ such that any vector of norm $<\epsi$ is the midpoint of a unique \textup{(}at most one\textup{)} chord of the unit circle,
\item There is an $\epsi>0$ such that any vector of norm $<\epsi$ is the sum of two unit vectors in a unique \textup{(}at most one\textup{)} way,
\item There is an $\epsi>0$ such that any two unit circles with centres at distance $<\epsi$ intersect in exactly \textup{(}at most\textup{)} two points,
\item Any two circles intersect in at most two points,
\item Any three points are contained in at most one circle.
\end{enumerate}
\end{proposition}

\proof
Essentially the only nontrivial part that remains to be proven is that given any three points $T=\{\vx_1,\vx_2,\vx_3\}$, and any positive homothetic copy $T'=\{\vx_1',\vx_2',\vx_3'\}$ of them, the resulting six points cannot be in strictly convex position, i.e., one of them is in the convex hull of the remaining five.
We may assume without loss of generality that $T$ is not collinear.

We first consider the case where $T'$ is not a translate of $T$.
Let $\vo$ be the centre of homothety.
If $\vo\in\conv T$, then $T\cup T'$ is clearly not in strictly convex position.

\begin{figure}
\begin{minipage}[b]{5.5cm}
\begin{center}
\includegraphics{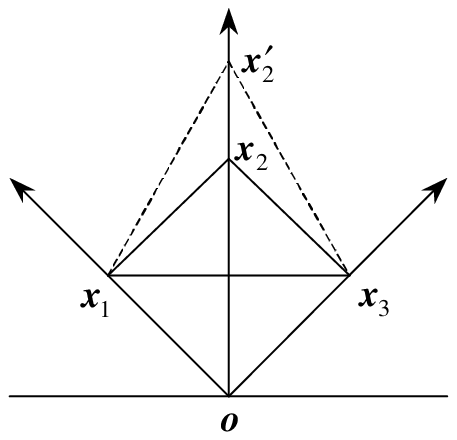}
\caption{$\vo$ and $\vx_2$ on opposite sides of $\Line{\vx_1\vx_3}$}\label{fig7}
\end{center}
\end{minipage}
\hfill
\begin{minipage}[b]{5.5cm}
\begin{center}
\includegraphics{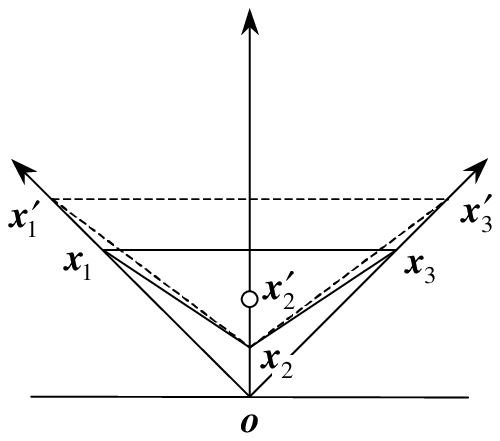}
\caption{$\vo$ and $\vx_2$ on the same side of $\Line{\vx_1\vx_3}$}\label{fig8}
\end{center}
\end{minipage}
\end{figure}

So assume that $\vo\not\in\conv T$.
Find a line $\ell$ through $\vo$ such that $T$ is in an open half-plane bounded by $\ell$.
See Figures~\ref{fig7} and \ref{fig8}.

We now show that some two rays $\Ray{\vo\vx_i}$ and $\Ray{\vo\vx_j}$ coincide.
If the rays $\Ray{\vo\vx_i}$ are all distinct, some ray is between the other two.
Assume without loss of generality that $\Ray{\vo\vx_2}$ is between $\Ray{\vo\vx_1}$ and $\Ray{\vo\vx_3}$.
Also assume without loss of generality that the homothety factor is $>1$.
If $\vo$ and $\vx_2$ are on opposite sides of the line $\Line{\vx_1\vx_3}$, then $\vx_2\in\interior\conv\{\vx_1,\vx_2',\vx_3\}$, a contradiction.
If $\vo$ and $\vx_2$ are on the same side of the line $\Line{\vx_1\vx_3}$, then $\vo$ and $\vx_2'$ are on the same side of $\Line{\vx_1'\vx_3'}$, and $\vx_2'\in\interior\conv\{\vx_1',\vx_2,\vx_3'\}$, a contradiction.

Thus some two rays $\Ray{\vo\vx_i}$ and $\Ray{\vo\vx_j}$ coincide, and then $\{\vx_i,\vx_j,\vx_i',\vx_j'\}$ is collinear.

The case where $T'$ is a translate of $T$ is similar:
If two of the lines $\Line{\vx_i\vx_i'}$ coincide, then some $\{\vx_i,\vx_j,\vx_i',\vx_j'\}$ is collinear.
Otherwise one of the lines is between the other two, say $\Line{\vx_2\vx_2'}$ between $\Line{\vx_1\vx_1'}$ and $\Line{\vx_3\vx_3'}$.
If $\vx_1,\vx_2,\vx_3$ are not collinear, then we obtain a contradiction as before by considering whether $\vx_2$ and $\Line{\vx_1'\vx_3'}$ are on opposite sides of the line $\Line{\vx_1\vx_3}$ or not.
\qed

Holub's second characterization is in terms of the \emph{bisector} of two points $\vx$ and $\vy$:
\[ B(\vx,\vy) := \{\vz : \length{\vx\vz}=\length{\vy\vz}\}. \]
As Holub did not give a proof in \cite{MR52:1263} and seemingly neither in a later paper, we include a proof.
We first give a related characterization for Minkowski spaces.

\begin{proposition}\label{bisector1}
A Minkowski space is strictly convex iff for all distinct points $\vx$ and $\vy$, and for all lines $\ell$ parallel to $\vy-\vx$, $\ell$ intersects $B(\vx,\vy)$ in exactly \textup{(}at most\textup{)} one point.
\end{proposition}

\proof
Suppose that $\vz_1$ and $\vz_2$ are two points in $B(\vx,\vy)$ with $\Line{\vz_1\vz_2}$ parallel to $\Line{\vx\vy}$.
Then the sum of the diagonals of the convex quadrilateral $\vx\vy\vz_1\vz_2$ equals the sum of the two opposite sides, and by Proposition~\ref{quadlemma} the space is not strictly convex.

Conversely, suppose the space is not strictly convex.
Let $\Segment{\va\vb}$ be a segment on the unit circle.
Let $\vx=\va, \vy=\frac{1}{2}(\va+\vb), \vz_1=\vo$ and $\vz_2=\frac{1}{2}(\va-\vb)$.
Then $\vz_1,\vz_2\in B(\vx,\vy)$, but $\Line{\vz_1\vz_2}$ is parallel to $\Line{\vx\vy}$.
\qed

A local version of the previous proposition is

\begin{corollary}
Given distinct $\vx,\vy$ in a Minkowski space, then all lines parallel to $\Line{\vx\vy}$ intersect $B(\vx,\vy)$ in at most one point iff the unit sphere does not contain a segment parallel to $\Line{\vx\vy}$.
\end{corollary}

\begin{proposition}[Holub \cite{MR52:1263}]\label{holub}
A Minkowski plane is strictly convex iff for any distinct $\vx,\vy$ and any $\vz\in B(\vx,\vy)$ the bisector of $\vx$ and $\vy$ is contained in the double cone of $\vx$ and $\vy$ with apex $\vz$, i.e.,
\[ B(\vx,\vy) \subseteq \{\vz+\lambda(\vx-\vz)+\mu(\vy-\vz) : \lambda\mu\geq 0\}. \]
\end{proposition}

\begin{figure}
\begin{center}
\includegraphics{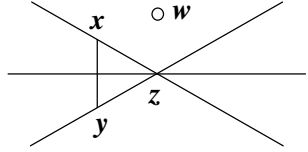}
\end{center}
\caption{Proof of Proposition~\ref{holub}}\label{fig9}
\end{figure}

\proof
Suppose that there is $\vw\in B(\vx,\vy)$ outside the double cone, say
$\vw=\vz+\lambda(\vx-\vz)+\mu(\vy-\vz)$ with $\lambda>0, \mu<0$ (see Figure~\ref{fig9}).
Using the convex quadrilateral $\vx\vy\vz\vw$ we obtain that the plane is not strictly convex as in the proof of Proposition~\ref{bisector1}.

For the converse the same example may be used as in the proof of Proposition~\ref{bisector1}.
\qed

We may again formulate a local version:

\begin{corollary}
Given distinct $\vx,\vy$ in a Minkowski plane, then $B(\vx,\vy)$ is contained in all double cones with apex $\vz$ iff there is no segment on the unit circle parallel to $\Line{\vx\vy}$.
In this case we in fact have that $B(\vx,\vy)$ equals the intersection of all the double cones as $\vz$ ranges over $B(\vx,\vy)$.
\end{corollary}

We now discuss a measure of non-strict convexity of Minkowski planes, introduced by Brass \cite{MR97c:52036}.
Let $\lambda(M)$ be the length of the longest line segment contained in the boundary of the unit circle of the Minkowski plane $M$.

\begin{proposition}\label{lambda}
$0\leq\lambda(M)\leq 2$, with equality on the left side iff $M$ is strictly convex, and equality on the right side iff $M$ is rectilinear, in which case the segments of length two ``fill up'' the unit circle.
\end{proposition}

This observation has been made more than once; see e.g.\ Brass \cite{MR97c:52036}.
Thus in the rectilinear plane there are two pairs of segments of length two on the unit circle.
Brass has shown that there are at most two pairs of segments of length $>1$.
Here we show the following sharpening:
\begin{proposition}\label{longsegmentprop}
On the unit circle of a Minkowski plane there are at most three pairs of segments of length at least $1$.
If there are three pairs of segments of length at least $1$, then the unit disc must a hexagon with vertices $\pm\vx_1,\pm\vx_2,\pm\lambda(\vx_1+\vx_2)$ for some $\lambda\in(\frac{1}{2},1]$, and at least two pairs are of length exactly $1$.
\end{proposition}

The proof will be given after we discuss the useful concept of conjugate diameters (Section~\ref{conjdiamsection}).

\subsection{Interlude: Intersection of homothets of a fixed convex curve}\label{intercircsect}
Although this section is not strictly Minkowski geometry, but rather convex geometry, the results are important also in Minkowski geometry; we need them for example in the proof of Proposition~\ref{innermetric} below.

As discussed in Section~\ref{scsect}, two circles in a Minkowski plane intersect in at most two points iff the plane is strictly convex.
This remains true for homothets of a closed convex curve that is not centrally symmetric, with the same proof (as noted in \cite{Mayer} and \cite{Buter}).
We now examine the case when the plane is not necessarily strictly convex.
In general, the intersection of two circles is always the union of two segments, which are either disjoint or intersect in a common endpoint, where a segment may degenerate to a point or the empty set.
This was proved by Gr\"unbaum \cite{MR21:2209} and later also by Banasiak \cite{MR89k:46021}.
Again, this result remains true if we do not assume central symmetry, as shown by Shiffman \cite{MR14:632d} (as mentioned in \cite{Novikoff}), and also by Sch\"affer \cite[4B]{MR57:7120}.
In the following proposition is also included a statement on where the different pieces of the homothets lie relative to each other; this generalizes a lemma of Sch\"affer \cite[Lemma 4.3]{MR36:1959}.

\begin{proposition}\label{intersectprop}
Let $C$ be a compact convex disc with boundary the closed convex curve $\gamma$, and $C'$ be a positive homothet \textup{(}which may be a translate\textup{)} with boundary $\gamma'$.
Then $\gamma\cap\gamma'$ is a union of two segments, each of which may degenerate to a point or the empty set.

Suppose furthermore that this intersection consists of two connected non-empty components $A_1,A_2$.
Then the centre of homothety is outside $C\cup C'$, and the two lines of homothety supporting $C\cap C'$ intersect $C\cap C'$ in exactly $A_1$ and $A_2$.

Choose a point $\vp_i$ from each component $A_i$ and let $T:C\to C'$ be the positive homothety mapping $C$ onto $C'$.
Let $\vc_i=T^{-1}\vp_i$ and $\vc_i'=T\vp_i$ for $i=1,2$.
Let $\gamma_1$ \textup{(}$\gamma_2$\textup{)} be the part of $\gamma$ on the same side \textup{(}opposite side\textup{)} of $\Line{\vp_1\vp_2}$ as $\vc_1$ and $\vc_2$; similarly for $\gamma'$.
See Figure~\ref{fig10}.

\begin{figure}
\begin{center}
\includegraphics{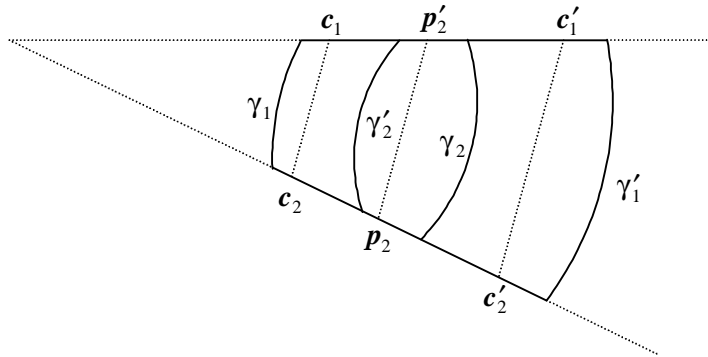}
\end{center}
\caption{Intersecting two homothets of a convex curve}\label{fig10}
\end{figure}

Then $\gamma_2\subseteq\conv \gamma_1'$ and $\gamma_2'\subseteq\conv \gamma_1$.
\end{proposition}

We reiterate for translates:

\begin{proposition}
Let $C$ be a compact convex disc with boundary the closed curve $\gamma$, and $C+\vv$ be a translate of $C$ with boundary $\gamma'$.
Then $\gamma\cap\gamma'$ is a union of two segments parallel to the direction of translation, each of which may degenerate to a point or the empty set.

Suppose that this intersection consists of two connected non-empty components $A_1,A_2$.
Then the two lines of translation supporting $C\cap C'$ intersect $C\cap C'$ in exactly $A_1$ and $A_2$.

Choose a point $\vp_i$ from each component $A_i$ and let $\vc_i=\vp_i-\vv$ and $\vc_i'=\vp_i+\vv$ for $i=1,2$.
Let $\gamma_1$ \textup{(}$\gamma_2$\textup{)} be the part of $\gamma$ on the same side \textup{(}opposite side\textup{)} of $\Line{\vp_1\vp_2}$ as $\vc_1$ and $\vc_2$; similarly for $\gamma'$.
Then $\gamma_2\subseteq\conv \gamma_1'$ and $\gamma_2'\subseteq\conv \gamma_1$.
\end{proposition}

We note that it can also be proved that the intersection of any number of homothets of a fixed closed convex curve has at most two components, except if the curve is a triangle, in which case it is also possible for the intersection of at least three homothets to be three points homothetic to the vertices of the original triangle.

\begin{figure}
\begin{center}
\includegraphics{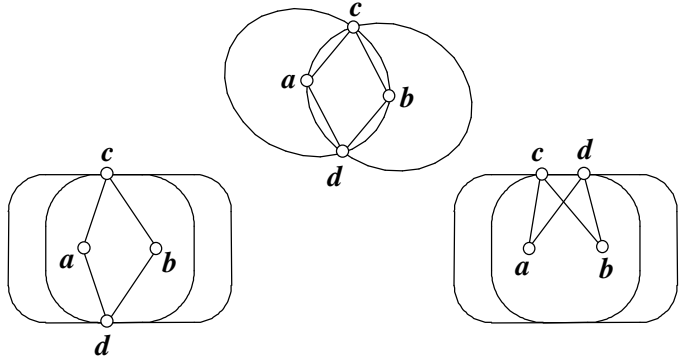}
\end{center}
\caption{The unit distance graph $K_{2,2}$}\label{brassfig}
\end{figure}

From the result on the intersection of circles it can be seen exactly when a four-cycle can occur as a unit distance graph, which is used to analyse unit distance graphs \cite{MR97c:52036}.
See Figure~\ref{brassfig}.

\begin{lemma}\label{fourunit}
Consider a quadrilateral with vertices $\va,\vb,\vc,\vd$ \textup{(}not necessarily in this order\textup{)} and sides of unit length.
Then the vertices may be relabled such that one of the following three cases occur:
\begin{enumerate}
\item $\va\vc\vb\vd$ is a parallelogram.
\item There is a segment on the unit circle parallel to $\Segment{\vc\vd}$, $\length{\va\vb}=2$ and $\length{\va\vc}=\length{\va\vd}=\length{\vb\vc}=\length{\vb\vd}=1$; also $\length{\vc\vd}\leq 2$ with equality implying that the plane is rectilinear.
\item There is a segment of length at least $\length{\va\vb}+\length{\vc\vd}$ on the unit circle parallel to $\Segment{\va\vb}$ and $\Segment{\vc\vd}$, and $\length{\va\vb}<2$ and $\length{\vc\vd}<2$.
\end{enumerate}
\end{lemma}

It is easily seen by Proposition~\ref{intersectprop} that the set-theoretic difference of two homothets of a convex disc is always connected.
This is generalized to higher dimensions by considering two-dimensional sections of the difference of two homothets of a convex body through two points that need to be connected (Banasiak \cite{MR89k:46021}).
However, the intersection of two spheres is complicated: if it is assumed that two unit spheres always have an intersection contained in a hyperplane, then the space must be Euclidean (Goodey \cite{MR84e:52006}).
A higher dimensional analogue of Proposition~\ref{intersectprop} is that the intersection of two spheres (or homothets of the boundary of a convex body) always consists of at most one connected component (announced by Novikoff \cite{Novikoff}).

A converse of the fact that the intersection of two homothets of a convex curve has at most two connected components is the following result:
Let $K_1$ and $K_2$ be two planar convex bodies, such that for any translate $K_2'$ of $K_2$ with $K_2'\neq K_1$ and the property that the interiors of $K_1$ and $K_2'$ intersect, we have that the boundaries of $K_1$ and $K_2'$ intersect in exactly two connected components.
Then $K_1$ and $K_2$ are translates (Goodey and Woodcock \cite{MR84e:52014}).
For related results see also \cite{MR84e:52007}.

\subsection{Generalization of the triangle inequality to convex curves}
We now give an elementary proof that in any Minkowski plane, if one convex curve is contained in another, it is not longer than the containing curve.
Archimedes in his work Sphere and Cylinder, Book I \cite{MR88i:01164}, took this statement as a postulate (together with his now famous axiom that any two magnitudes are comparable using integer multiples).
Perhaps the most elegant way of proving this is by mixed volumes (as in Thompson \cite[Remark on p.~121]{MR97f:52001}) or by Crofton formulae \cite{AlvarezDuran}, but it is not difficult to give an elementary treatment, and to characterize equality, which delivers a small surprise (Corollary~\ref{strictineq}).
Gr\"unbaum \cite{MR34:1925} also uses this inequality for Minkowski planes.

\begin{lemma}\label{pathineq}
Let $\va_1\dots\va_n$ be a convex polygonal path contained in the triangle $\va_1\vb\va_n$.
Then $\sum_{i=1}^n\length{\va_i\va_{i+1}}\leq\length{\va_1\vb}+\length{\vb\va_n}$.
\end{lemma}

\begin{figure}
\begin{center}
\includegraphics{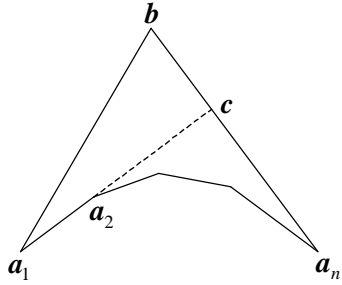}
\end{center}
\caption{Proof of Lemma~\ref{pathineq}}\label{archimedes}
\end{figure}

\proof
Induction on $n\geq 2$.
The case $n=2$ is the triangle inequality.
In the induction case $n\geq 3$, let $\Ray{\va_1\va_2}$ intersect the segment $\Segment{\vb\va_n}$ in $\vc$.
See Figure~\ref{archimedes}.

Then
\begin{eqnarray*}
\sum_{i=1}^{n-1}\length{\va_i\va_{i+1}} & = & \length{\va_1\va_2}+\sum_{i=2}^{n-1}\length{\va_i\va_{i+1}} \\
& \leq & \length{\va_1\va_2}+\length{\va_2\vc}+\length{\vc\va_n} \quad \mbox{(Induction)} \\
& \leq & \length{\va_1\vb}+\length{\vb\vc}+\length{\vc\va_n} \quad \mbox{(Triangle inequality)}.
\end{eqnarray*}
\qed

\begin{lemma}\label{curvetrineq}
Let $\gamma$ be a convex curve from $\va$ to $\vc$ contained in triangle $\va\vb\vc$.
Then $\length{\gamma}\leq\length{\va\vb}+\length{\vb\vc}$.
\end{lemma}

\proof
Approximate $\gamma$ by inscribed convex polygons (which are still contained in the triangle by convexity), use Lemma~\ref{pathineq}, and take the supremum.
\qed

\begin{lemma}\label{almostconvexcurveineq}
Let $\gamma$ be a convex curve from $\va_1$ to $\va_n$, contained in a convex polygon $\va_1\dots\va_n$.
Then $\length{\gamma}\leq\sum_{i=1}^{n-1}\length{\va_i\va_{i+1}}$.
\end{lemma}

\proof
We use induction on $n\geq 3$.
The case $n=3$ is Lemma~\ref{curvetrineq}.
For $n>3$, take an edge of the polygon not containing $\va_1$ or
$\va_n$ and translate it until it supports $\curve{\gamma}$ at $\vb$, say.
Then we have two polygons $\va_1\dots\vb$ and $\vb\dots\va_n$, and we
may use the induction hypothesis and the triangle inequality.
\qed

\begin{theorem}\label{convexcurveineq}
Let $\gamma_1$ and $\gamma_2$ be convex curves from $\va$ to
$\vb$ such that $\conv\gamma_1\subseteq\conv\gamma_2$.
Then $\length{\gamma_1}\leq\length{\gamma_2}$, with equality iff there
exists a \textup{(}possibly infinite\textup{)} sequence $\va_1,\va_2,\dots$ common to
$\gamma_1$ and $\gamma_2$, and such that the arcs on
$\gamma_1$ and $\gamma_2$ from $\va_i$ to $\va_{i+1}$ are both
metric segments.
\end{theorem}

\begin{figure}
\begin{center}
\includegraphics{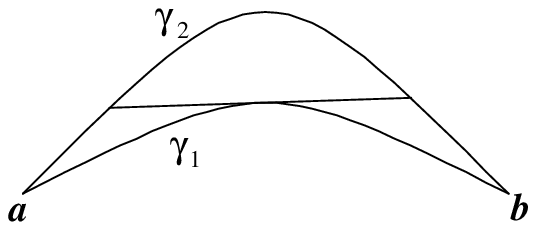}
\end{center}
\caption{$\length{\gamma_1}\leq\length{\gamma_2}$}\label{archimedes2}
\end{figure}

\proof
To obtain the inequality, approximate $\gamma_2$ by circumscribed
polygons, use Lemma~\ref{almostconvexcurveineq} and take the limit.

We now assume that $\length{\gamma_1}=\length{\gamma_2}$.
The intersection $\gamma_1\cap\gamma_2$ is a closed subset of a
(topological) interval (or circle).
By considering the maximal open intervals contained in
$\gamma_1\cap\gamma_2$ as well as its complement, we obtain a
partition of $\gamma_1$ and $\gamma_2$ into corresponding pieces
that are either identical or coincide only at endpoints.
We may therefore assume without loss of generality that $\curve{\gamma_1}$ and
$\gamma_2$ coincide only at endpoints.
If we cut off a piece of $\gamma_2$ using a supporting line of
$\gamma_1$, we obtain a metric segment (from the inequality $\length{\gamma_1}\leq\length{\gamma_2}$).
See Figure~\ref{archimedes2}.
By Proposition~\ref{msegchar} all chords of this metric segment are
contained in a segment of $S$.
Using all supporting lines of $\gamma_1$, we obtain that all
chords of $\gamma_2$ are contained in the same edge of $S$.
Again by Proposition~\ref{msegchar}, $\gamma_2$ is a metric segment,
hence $\gamma_1$ is also a metric segment.
\qed

\begin{corollary}\label{strictineq}
Let $\gamma_1$ and $\gamma_2$ be convex curves from $\va$ to $\vb$ such that $\conv\gamma_1\subseteq\conv\gamma_2$, $\length{\va\vb}<\length{\gamma_1}$, and $\gamma_1$ and $\gamma_2$ only have endpoints in common.
Then $\length{\gamma_1}<\length{\gamma_2}$.
\end{corollary}

It is now simple to prove the following lemma used by Sch\"affer \cite[4E]{MR57:7120}.

\begin{proposition}\label{innermetric}
Let $\vp,\vp',\vq',\vq$ be points in the given order on a convex curve $\gamma$ such that $\Line{\vp\vq}$ and $\Line{\vp'\vq'}$ are parallel.
Let $\gamma_1$ be the curve from $\vp$ to $\vq$ containing $\vp'$ and $\vq'$, and $\gamma_2$ the part of $\gamma_1$ from $\vp'$ to $\vq'$.
Then \[ \frac{\length{\gamma_1}}{\length{\vp\vq}} \geq \frac{\length{\gamma_2}}{\length{\vp'\vq'}}. \]
\end{proposition}

\begin{figure}
\begin{center}
\includegraphics{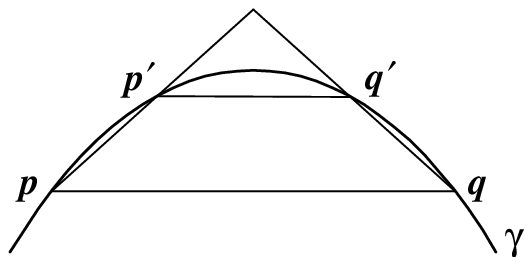}
\end{center}
\caption{Proof of Proposition~\ref{innermetric}}\label{schaeffer}
\end{figure}

\proof
Using Proposition~\ref{intersectprop}, the image $\gamma_1'$ of $\gamma_1$ under the positive homothety mapping $\Segment{\vp\vq}$ to $\Segment{\vp'\vq'}$ contains $\gamma_2$.
See Figure~\ref{schaeffer}.
By Theorem~\ref{convexcurveineq} it follows that
\[ \length{\gamma_2}\leq\length{\gamma_1'} = \length{\gamma_1}\cdot\frac{\length{\vp'\vq'}}{\length{\vp\vq}}.\]
\qed

Sch\"affer \cite{MR36:1959} studied the \emph{inner metric} of the unit sphere $S$ of a Minkowski space:
For unit vectors $\vp,\vq$, $\delta(\vp,\vq)$ is the infimum (minimum in finite dimensional spaces) of the lengths of all curves on $S$ joining $\vp$ and $\vq$.
He proved that the inner metric and the metric induced by the norm are equivalent using the following planar inequality \cite[Theorem~4.4]{MR36:1959}, which is an immediate corollary of the previous propostion.

\begin{proposition}
In a Minkowski plane $M$, for any unit vectors $\vp',\vq'$,
\[\delta(\vp',\vq')\leq\frac{1}{2}\Pi(M)\length{\vp'\vq'}. \]
\end{proposition}
\proof
Let $\Segment{\vp\vq}$ be the diameter of the unit circle parallel to $\Segment{\vp'\vq'}$.
Note that $\delta(\vp',\vq')$ is the length of the arc of the unit circle from $\vp'$ to $\vq'$ on the side of $\Segment{\vp'\vq'}$ opposite $\Segment{\vp\vq}$.
Now apply Proposition~\ref{innermetric}.
\qed

\subsection{The monotonicity lemma}\label{monotonsubsect}
For a fixed point $\vp$ on the unit circle and a variable point $\vx$, the length $\length{\vp\vx}$ is non-decreasing as $\vx$ moves on the unit circle from $\vp$ to $-\vp$.
See Figure{monfig}.
\begin{figure}
\begin{center}
\includegraphics{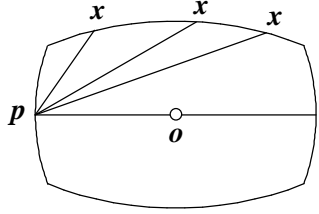}\label{monfig}
\end{center}
\caption{Monotonicity Lemma}
\end{figure}
This has often been assumed as intuitively obvious and has sometimes been given with complicated proofs (Sch\"affer \cite{MR57:7120}, Thompson \cite[Lemma 4.1.2]{MR97f:52001}, Doyle-Lagarias-Randall \cite{MR93e:52038}).
Thompson writes ``a more straightforward proof would be preferable\dots''.
Gr\"unbaum \cite{MR21:2209} gives a proof of the following strengthening: if $\va,\vb,\vc,\vd$ are on the unit circle in this order, then $\length{\vb\vc}<\length{\va\vd}$, unless $\length{\vb\vc}=\length{\va\vd}=2$ in which case the plane is not strictly convex, and the points are on opposite pairs of line segments in the unit circle.

Here we give a simple proof of a generalization ($\vp$ does not have to be on the unit circle), mentioned without proof by Alonso and Ben{\'\i}tez \cite[Lemma 1]{MR90k:46047}, and characterize equality.
The proof is a generalization of a proof of Valentine \cite{MR84m:46023} that two unit circles intersect in at most two points.

\begin{proposition}[Monotonicity lemma]\label{monlemma}
Let $\va,\vb\,\vc\neq\vo$, $\va\neq\vc$, with $\Ray{\vo\vb}$ between $\Ray{\vo\va}$ and $\Ray{\vo\vc}$, and suppose that $\length{\vo\vb}=\length{\vo\vc}$.
Then $\length{\va\vb}\leq\length{\va\vc}$, with equality iff either
\begin{enumerate}
\item $\vb=\vc$,
\item or $\vo$ and $\vb$ are on opposite sides of $\Line{\va\vc}$, and $\Segment{\un{\vc-\va}\;\;\un{\vb}}$ is a segment on the unit circle,
\item or $\vo$ and $\vb$ are on the same side of $\Line{\va\vc}$, and $\Segment{\un{\vc-\va}\;\;\un{-\vc}}$ is a segment on the unit circle.
\end{enumerate}
In particular, if the plane is strictly convex, then we always have strict inequality.
\end{proposition}

\begin{figure}
\begin{minipage}[b]{5.5cm}
\begin{center}
\includegraphics{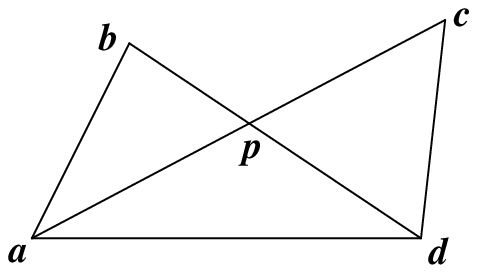}
\end{center}
\caption{Proof of Monotonicity Lemma}\label{fig11a}
\end{minipage}
\hfill
\begin{minipage}[b]{5.5cm}
\begin{center}
\includegraphics{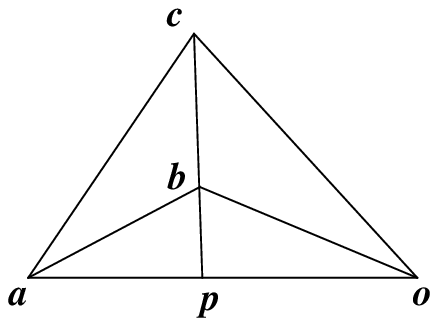}
\caption{Proof of Monotonicity Lemma}\label{fig11b}
\end{center}
\end{minipage}
\end{figure}

\proof
If $\vo$ and $\vb$ are on opposite sides of $\Line{\va\vc}$, let $\Segment{\vo\vb}$ and $\Segment{\va\vc}$ intersect in $\vp$.
See Figure~\ref{fig11a}.

Then
\begin{eqnarray*}
\length{\vo\vb}+\length{\va\vc} & = & (\length{\va\vp} + \length{\vp\vb}) + (\length{\vo\vp} + \length{\vp\vc}) \\
& \geq & \length{\va\vb} + \length{\vo\vc} \quad \mbox{(Triangle inequality).}
\end{eqnarray*}
The case of equality may now be analysed using Proposition~\ref{trieq}.

If $\vo$ and $\vb$ are on the same side of $\Line{\va\vc}$, let $\Ray{\vc\vb}$ intersect $\Segment{\vo\va}$ in $\vp$.
See Figure~\ref{fig11b}.
Since $\length{\vo\vb}=\length{\vo\vc}$, Lemma~\ref{maxineq} gives that $\length{\vo\vp}\geq\length{\vo\vc}$.
Thus $\length{\va\vc}+\length{\vo\vp}\geq\length{\va\vc}+\length{\vo\vc}\geq\length{\va\vo}=\length{\va\vp}+\length{\vo\vp}$ by the triangle inequality, hence $\length{\va\vc}\geq\length{\va\vp}$, and $\length{\va\vc}\geq\length{\va\vb}$ by Lemma~\ref{maxineq}.
Again, equality can be analysed using Proposition~\ref{trieq}.
\qed

\begin{figure}[b]
\begin{center}
\includegraphics{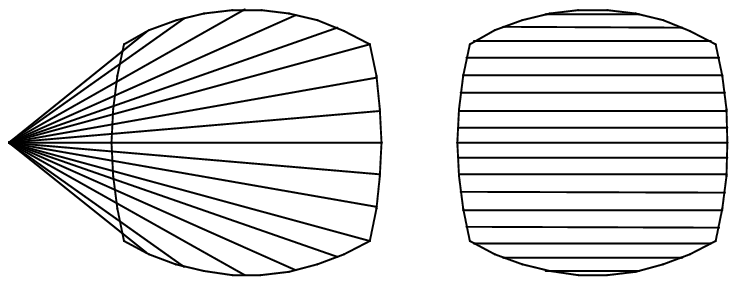}
\end{center}
\caption{}\label{fig18}
\end{figure}

\begin{figure}
\begin{center}
\includegraphics{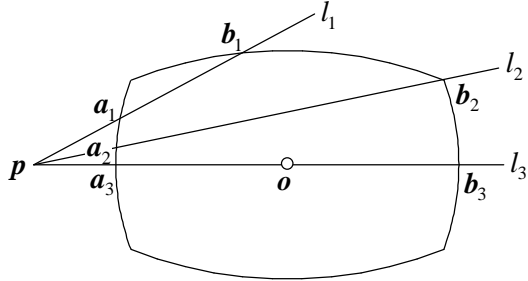}
\end{center}
\caption{Proof of Corollary~\ref{unimod}}\label{fig19}
\end{figure}

As a simple corollary we have that a pencil of lines with point of concurrency outside the unit circle intersects the unit circle in segments of which the lengths are unimodal.
See Figures~\ref{fig18} and \ref{fig19}.

\begin{corollary}\label{unimod}
Let $\norm{\vp}>1$.
Let $\ell_1,\ell_2,\ell_3$ be lines through $\vp$ with $\ell_2$ between $\ell_1$ and $\ell_3$, and $\ell_3$ a diameter of the unit circle.
Let $\ell_i$ intersect the unit circle in $\va_1$ and $\vb_i$, $i=1,2,3$.
Then $\length{\va_1\vb_2}\leq\length{\va_2\vb_2}$.
\end{corollary}

Also, by specializing the above proof of the Monotonicity Lemma to strictly convex planes, we obtain Valentine's proof \cite{MR84m:46023} that two unit circles intersect in at most two points (we also do not need to assume that the circles have equal radii).
Gr\"unbaum and Kelly \cite[Theorem 3]{MR39:6180}) use a similar method as the first case of the above proof to obtain a monotonicity result for curves of \emph{constant diameter}, i.e., each point is an endpoint of a diameter, in strictly convex Minkowski planes.
It would be of interest to characterize the curves for which there is a monotonicity result; perhaps the curves of constant diameter are the only such ones.
Heppes \cite{MR22:1846} (see \cite{MR39:6180}) has shown that in the Euclidean plane, if a monotonicity result holds, then the curve has constant width (which in the Euclidean case is equivalent to constant diameter).

\section{Equilateral triangles and affine regular hexagons}
\begin{proposition}\label{eqtri}
Given any segment $\Segment{\vp\vq}$ in a Minkowski plane, and a half-plane bounded by the line $\Line{\vp\vq}$, there exists a point $\vr$ in the half plane such that $\triangle\vp\vq\vr$ is an equilateral triangle.
The point $\vr$ is unique iff the longest segment in the unit circle parallel to $\Line{\vp\vq}$ has length at most $1$.
\end{proposition}

\begin{figure}[b]
\begin{center}
\includegraphics{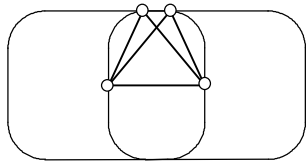}
\end{center}
\caption{Non-uniqueness in the construction of an equilateral triangle}\label{fig20}
\end{figure}

\proof
The proof is exactly the same as the proof for the Euclidean plane going back to Euclid:
consider the two circles which have one endpoint of the segment as a centre, and the other endpoint on the boundary.
Any intersection point $\vr$ of the two circles gives the third point of the equilateral triangle.
If two points $\vr_1$ and $\vr_2$ are both intersection points of the two circles in the chosen half plane, then by Lemma~\ref{fourunit} we obtain a segment of length $1+\length{\vr_1\vr_2}/\length{\vp\vq}$ on the unit circle, parallel to $\Segment{\vp\vq}$. 
See Figure~\ref{fig20}.
The converse is clear.
\qed

Thus we get uniqueness in a larger class of Minkowski planes than the strictly convex ones, namely those in which the measure of non-strict convexity $\lambda(M)\leq 1$.
Also, by Proposition~\ref{longsegmentprop} it follows that it is in at most two directions of the given segment that we get non-uniqueness of the equilateral triangle, for any Minkowski plane.

From the construction of an equilateral triangle it is possible to construct a hexagon inscribed in the unit circle and with unit side lengths, exactly as in Euclidean geometry.
This construction was noticed very early (Go{\l}ab \cite{Golab}, P. J. Kelly \cite{MR12:525i}).

\begin{figure}
\begin{center}
\includegraphics{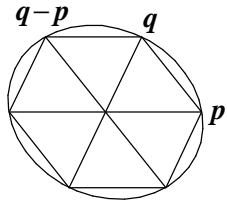}
\end{center}
\caption{An affine regular hexagon inscribed to the unit circle}\label{fig21}
\end{figure}

\begin{proposition}
Let $\triangle\vo\vp\vq$ be an equilateral triangle.
Then the hexagon with vertices $\pm\vp,\pm\vq,\pm(\vp-\vq)$ is an affine image of a Euclidean equilateral hexagon, and all sides have the same length.
\end{proposition}
See Figure~\ref{fig21}.
If $\vp_1,\dots,\vp_6$ are the six vertices in order of some inscribed regular hexagon, then we call $\vo\vp_i$ a \emph{radius}, $\vp_i\vp_{i+3}$ a \emph{diameter}, and $\vp_i\vp_{i+2}$ a \emph{diagonal}.
Given a direction for a diameter or radius of such a hexagon, we have the same discussion of uniqueness as for the equilateral triangle:
there is an essentially unique hexagon iff the unit ball does not contain a segment of length $>1$ parallel to the given direction.
Lassak \cite{MR91a:52001} proves that if the direction of a diagonal is given, there is always a unique hexagon.

For non-strictly convex curves one can also always find an inscribed affine regular hexagon, but the direction may not be prescribed anymore (discovered independently by various authors; see Gr\"unbaum \cite[p.\ 242--243]{MR27:6187}).

The following two finer results have been used as intermediate steps to various results in planar Minkowski geometry.

Ohmann \cite{MR14:76a} shows that any convex disc in a Minkowski plane has a symmetric circumscribed hexagon whose sides are parallel to some hexagon circumscribed to the unit ball $B$ and supporting the unit ball at the vertices of an inscribed affine regular hexagon.

Lassak \cite{MR91a:52001} also proves that the unit circle has an inscribed affine regular hexagon $\vp_1\dots\vp_6$ such that the line through $\vp_i$ parallel to $\Line{\vp_{i-1}\vp_{i+1}}$ supports the unit disc.

In three-dimensional Minkowski spaces the following results are known.
Sch\"affer \cite{MR38:1610} shows with topological methods that a centrally symmetric simple closed curve on the unit sphere of a three-dimensional Minkowski space has a concentric inscribed affine regular hexagon.
Laugwitz \cite{MR94c:51030} shows that for any two unit vectors $\vu_1,\vu_2$ in a three-dimensional Minkowski space there exists a unit vector $\vu_3$ and affine regular hexagons $H_1, H_2$ inscribed in the unit ball such that $\vu_1,\vu_3\in H_1$ and $\vu_2,\vu_3\in H_2$.
Again the proof needs some topology.

See Ceder \cite{MR29:6391}, Ceder and Gr\"unbaum \cite{MR35:3544} for further results on inscribed (and circumscribed) hexagons.

We now discuss some applications of the hexagon construction.
First of all, it immediately gives the lower bound of $6$ to the circumference of the unit circle (see Section~\ref{circumferencesect}), as well as the lower bound of $6$ for the Hadwiger number of the unit disc (see Gr\"unbaum \cite{MR25:1492} for a discussion of the upper bound).
Chakerian \cite{MR34:4986} gives a discussion of how it is applied to prove that in the Minkowski plane a Reuleaux triangle is a body of fixed constant width of least area.
Lassak \cite{MR91a:52001} uses his strengthened hexagon construction mentioned above to obtain a two-dimensional Minkowski analogue of John's theorem on the approximation of convex bodies by the Euclidean ball.

We now give a proof of a special case of an inequality of Nordlander \cite{MR25:4329} using the hexagon construction.
Nordlander's inequality in general is proved by a continuous averaging argument.

\begin{proposition}
In any Minkowski plane,
\begin{equation}
\inf\{\norm{\vx+\vy} : \norm{\vx},\norm{\vy}\geq 1, \norm{\vx-\vy}\leq 1\} \leq \sqrt{3},\label{hexone}
\end{equation}
and
\begin{equation}
\sup\{\norm{\vx+\vy} : \norm{\vx},\norm{\vy}\leq 1, \norm{\vx-\vy}\geq 1\} \geq \sqrt{3}.\label{hextwo}
\end{equation}
\end{proposition}

\proof
We only prove inequality \eqref{hexone}, as \eqref{hextwo} is similar.
Let $\vx_i$, $i\in\Z_6$, be an affine regular hexagon inscribed in the unit circle, i.e., 
$\norm{\vx_0}=\norm{\vx_1}=\norm{\vx_1-\vx_0}=1$, and $\vx_{i+1}=\vx_i-\vx_{i-1}$.
Let $\vy_i=\frac{1}{\sqrt{3}}(\vx_i+\vx_{i+1})$ for all $i\in\Z_6$.
Then $\vy_{i+1}=\vy_i-\vy_{i-1}$ and $\vx_i=\frac{1}{\sqrt{3}}(\vy_i-\vy_{i-1})$ for all $i\in\Z_6$.

If $\norm{\vx_i+\vx_{i+1}}\leq\sqrt{3}$ for some $i$, then \eqref{hexone} is satisfied.
Otherwise $\norm{\vy_i}>1$ for all $i$.
Consider a $\vy_i$ of smallest norm.
Let $\vx=\frac{1}{\norm{\vy_i}}\vy_{i-1}$ and $\vy=\frac{1}{\norm{\vy_i}}\vy_{i-2}$.
Then $\norm{\vx},\norm{\vy}\geq 1$, $\norm{\vx-\vy}=\norm{\vy_{i-1}-\vy_{i-2}}/\norm{\vy_i}=1$, and
$\norm{\vx+\vy}=\norm{\vy_{i-1}+\vy_{i-2}}/\norm{\vy_i}=\sqrt{3}/\norm{\vy_i} < \sqrt{3}$,
proving \eqref{hexone}.
\qed

We now show how Lemma~\ref{maxineq} and the Monotonicity Lemma (Proposition~\ref{monlemma}) are utilized to sharpen the above result.

\begin{corollary}
In any Minkowski plane,
\begin{equation}
\inf\{\norm{\vx+\vy} : \norm{\vx}=\norm{\vy}=\norm{\vx-\vy}= 1\} \leq \sqrt{3},\label{hexthree}
\end{equation}
and
\begin{equation}
\sup\{\norm{\vx+\vy} : \norm{\vx}=\norm{\vy}=\norm{\vx-\vy}= 1\} \geq \sqrt{3}.\label{hexfour}
\end{equation}
\end{corollary}

\proof
Consider any $\vx,\vy$ such that $\norm{\vx},\norm{\vy}\geq 1$, $\norm{\vx-\vy}\leq 1$.
Let $d=\norm{\vx+\vy}$.
Assume without loss of generality that $\norm{\vx}\leq\norm{\vy}$.
Let $\vx_1=\un{\vx}$ and $\vy_1=\frac{1}{\norm{\vx}}\vy$.
Then $\norm{\vx_1-\vy_1}\leq 1$, $\norm{\vx_1}=1\leq\norm{\vy_1}$ and $\norm{\vx_1+\vy_1}\leq d$.
Let $\vy_2=\un{\vy_1}$.
By Lemma~\ref{maxineq}, $\norm{\vx_1-\vy_2}\leq\max\{\norm{\vx_1+\vy_1},\norm{\vx_1}\}\leq d$.
By continuity there exists a $\vy_3$ on the arc of the unit circle from $\vy_1$ to $-\vx_1$ such that $\norm{\vx_1-\vy_3}=1$.
By Proposition~\ref{monlemma}, $\norm{\vx_1+\vy_3}\leq\norm{\vx_1+\vy_2}\leq d$.

Inequality~\eqref{hexfour} is proved similarly.
\qed

Geometrically, the previous corollary says that there always exists a unit equilateral triangle with a median $\leq\sqrt{3}/2$ as well as one with a median $\geq\sqrt{3}/2$.
(A \emph{median} of a triangle is a segment from a vertex to the midpoint of the opposite edge.)
By continuity there then exists a unit equilateral triangle with a median of length exactly $\sqrt{3}/2$.
The Euclidean plane shows that the value of $\sqrt{3}/2$ is best possible.
Surprisingly, there are other Minkowski planes where the medians of all unit equilateral triangles are all $\sqrt{3}/2$, such as the plane with an affine regular hexagon as unit circle (L. M. Kelly \cite{MR53:6433}).
See Alonso and Ben{\'\i}tez \cite{MR90k:46047} for a discussion of equality in all other cases of Nordlander's inequality: for a countable set of cases, there is equality also for certain other Minkowski planes (such as those with regular polygons as unit circles), and for all other cases, there is a characterization of Euclidean space.
These results are further generalized by Alonso and Ull{\'a}n \cite{MR93m:46018}.

See Section~\ref{circumferencesect} for results on the ratio between the area of a unit equilateral triangle and the area of the unit disc.

\section{Equilateral sets}
\subsection{Four-point equilateral sets}
The following proposition is derived in many papers (e.g.\ Ra\u\i ko \cite{MR56:13114}, Chilakamarri \cite{MR92b:05036}, Brass \cite{MR97c:52036}).

\begin{proposition}
The maximum cardinality of an equilateral set is $4$ in the rectilinear plane, and $3$ in any other Minkowski plane.
\end{proposition}
\proof
We use the notation of the proof of Proposition~\ref{eqtri}.
In order to obtain an equilateral set of four points we have to find in the previous construction of an equilateral triangle two points $\vr_1$ and $\vr_2$ on the same side of $\Line{\vp\vq}$ such that $\length{\vr_1\vr_2}\geq \length{\vp\vq}$.
This means that on the unit circle we have to find a segment of length $1+\length{\vr_1\vr_2}/\length{\vp\vq} \geq 2$, i.e., we must have $\lambda(M_2)\geq 2$.
It follows from Proposition~\ref{lambda} that a four-point equilateral set is possible only in the rectilinear plane.
Also, since we now have $\length{\vr_1\vr_2}=\length{\vp\vq}$, a $5$-point equilateral set is impossible in the rectilinear plane.
\qed

The previous proposition can also be proved by a packing argument of Danzer and Gr\"unbaum \cite{MR25:1488}, by noticing that equilateral sets are antipodal sets (see Petty \cite{MR43:1051} and Soltan \cite{MR52:4127}).
This packing argument generalizes to higher dimensions (see below).

\subsection{Equilateral sets in higher dimensions}
Equilateral sets are considered by Baronti, Casini and Papini \cite{MR94m:46020} in relation to Chebyshev centres, medians, and barycentres.
Equilateral sets in Minkowski spaces are used to find area-minimizing surfaces, see Lawlor and Morgan \cite{MR95i:58051}.
See Section~\ref{circumferencesect} for results on the volume of equilateral simplices in higher dimensions.
We now survey the literature on the maximum cardinalities of these sets.

\subsubsection{Upper bounds}
Petty \cite{MR43:1051} and, independently, Soltan \cite{MR52:4127} show that equilateral sets are antipodal, and by a result of Danzer and Gr\"unbaum \cite{MR25:1488}, an antipodal set in $d$-dimensional vector space has at most $2^d$ elements, with equality iff the points form the vertex set of a parallelotope.
They also show that the existence of $2^d$ equilateral points implies that the unit ball is a parallelotope homothetic to the convex hull of the equilateral set.
This upper bound is also proved by F\"uredi, Lagarias and Morgan \cite{MR93d:52009} using the isodiametric inequality for Minkowski spaces due to Mel'nikov \cite{MR27:6191}.

Petty characterizes the Minkowski spaces in which equilateral sets are strictly antipodal; this class includes the strictly convex spaces.
Better upper bounds are known in some cases for strictly antipodal sets.
E.g., Gr\"unbaum \cite{MR28:2480} shows that in three-dimensional space a strictly antipodal set has at most five points.
It follows that in a strictly convex three-dimensional Minkowski space, an equilateral set can have at most five points.
There indeed exists such a Minkowski space, which is also smooth; see Lawlor and Morgan \cite[Example 3.4]{MR95i:58051}.
See also Morgan \cite{MR93h:53012}.

Soltan furthermore proves the following:
For any convex set $S$, let $b(S)$ be the least number of positive homothets of $S$ with homothety coefficient $<1$ that cover $S$.
Let $A$ be the convex hull of an equilateral set of cardinality $k$.
Then $k=b(A)$.
He also asks whether $b(B)$ equals the largest cardinality of an equilateral set.

It is easy to find $d+1$ equilateral points in $d$-dimensional $\ell_p$, and even $2d$ equilateral points in $d$-dimensional $\ell_1$.
However, it is not known whether these are upper bounds, except in three-dimensional $\ell_1$, where it is known that there are at most $6$ (Bandelt, Chepoi and Laurent \cite{MR99d:51017}).
See also \cite{Guy} where these problems are explicitly mentioned.

F\"uredi, Lagarias and Morgan \cite{MR93d:52009} conjecture that there exists an $\epsi>0$ such that for any strictly convex $d$-dimensional Minkowski space an equilateral set has cardinality at most $(2-\epsi)^d$.
They construct strictly convex spaces for each dimension such that there are equilateral sets of cardinality at least $(1.02)^d$.

In \cite{MR99k:52028} the following generalization of Petty's upper bound is conjectured (and the case $d=2$ is proved):
A subset of a $d$-dimensional Minkowski space in which at most $k$ distinct distances occur, has cardinality at most $(k+1)^d$, with equality only if the unit ball is a parallelotope.

\subsubsection{Lower bounds}
Petty showed that in a Minkowski space of at least three dimensions, there are always equilateral sets of size $4$.
He uses a topological result which essentially says that a punctured disc is not contractible.
Makeev \cite{Makeev} proves a stronger property in three dimensions, involving two norms.
Surprisingly, it is not known whether each $d$-dimensional Minkowski space $(d\geq 4)$ admits an equilateral set of $d+1$ points.
This question is asked by Ra\u\i ko \cite{MR56:13114}, Morgan \cite{MR93h:53012} and Thompson \cite{MR97f:52001}.
Also, Gr\"unbaum \cite[bottom of p.\ 242]{MR27:6187} asks for the existence of the difference body of a $d$-dimensional simplex inscribed in an arbitrary $d$-dimensional convex body.
For the case of symmetric convex bodies this is equivalent to the existence of $d+1$ equilateral points in the norm determined by the convex body.

The best result to date is that there is a function $f(d)$ tending to infinity such that any $d$-dimensional Minkowski space has an equilateral set of size $f(d)$ (Brass \cite{MR1720106}).
The proof uses Dvoretzky's theorem to reduce to an almost Euclidean space, and then Brouwer's theorem to turn an almost equilateral set into an equilateral set.

Petty \cite{MR43:1051} also gives an example of a $d$-dimensional space for each $d\geq 4$ where there is an equilateral set of $4$ points that is maximal, i.e., that cannot be extended to a larger equilateral set.
Thus a naive attempt at proving the above problem using induction seemingly does not work.

\section{Normality, Conjugate diameters, and Radon curves}\label{normalsection}
Normality, as defined in the Introduction, was introduced in a question of Carath\'eodory studied by Blaschke \cite{Blaschke} and Radon \cite{Radon, MR89i:01142a}, and introduced independently by Birkhoff \cite{Birkhoff}.
James \cite{MR9:42c} made a study of normality in normed spaces, linking it to strict convexity and smoothness (see Sections~\ref{scsect}).
From Proposition~\ref{uniquebody} it follows immediately that the normality relation determines the norm uniquely up to a positive constant.
This is also proved in detail by Sch\"opf \cite{MR98m:46018}.

The normality relation is not necessarily symmetric.
In fact, in Minkowski spaces of dimension at least $3$ normality is symmetric iff the space is Euclidean (the proof essentially going back to Blaschke \cite{Blaschke} and, independently, Birkhoff \cite{Birkhoff}; the final result, without any assumptions of strict convexity or smoothness, is due to James \cite{MR9:42d}; see Thompson \cite[\S3.2]{MR97f:52001} for a discussion.)

In two dimensions normality is symmetric iff the unit circle is a so-called Radon curve, introduced by Radon \cite{Radon, MR89i:01142a}.
These curves have many remarkable almost-Euclidean properties (see Section~\ref{radon} below).

James \cite{MR9:42c} showed that for any $\vx,\vy$ there exists an $a\in\R$ such that $\vx\normal a\vx+\vy$.
Also, if $\vx\normal a\vx+\vy$ and $\vy\normal b\vy+\vx$, then $\abs{ab}\leq 1$.
Necessity in the following characterization of symmetry of normality (in arbitrary dimensions) is due to James \cite{MR9:42c}.

\begin{proposition}
Normality in a Minkowski space is symmetric iff for any $\vx,\vy\neq\vo$ the following implication holds: if $\vx\normal a\vx+\vy$ and $\vy\normal b\vy+\vx$, then $ab\geq 0$.
\end{proposition}
Laugwitz gave a similar characterization in \cite{MR51:11069}.

Next we first discuss normality and conjugate diameters in Minkowski planes, where many of the essential ideas already occur.
Then we indicate higher-dimensional generalizations.

\subsection{Planes}
\subsubsection{Conjugate diameters}\label{conjdiamsection}
\begin{figure}
\begin{center}
\includegraphics{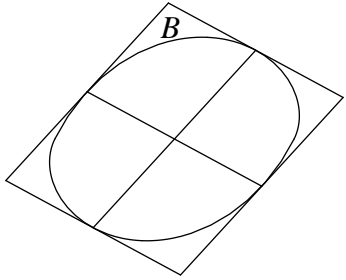}
\end{center}
\caption{A pair of conjugate diameters}\label{fig22}
\end{figure}

Two diameters of the unit circle of a Minkowski plane are \emph{conjugate} if their directions are mutually normal, i.e., $\vx\normal\vy$ and $\vy\normal\vx$, where $\vx$ and $\vy$ are the directions of the two diameters.
See Figure~\ref{fig22}.
It is not \emph{a priori} clear that conjugate diameters exist for any unit circle.
It seems that Auerbach \cite{Auerbach1, Auerbach2} (see also \cite[p.~238]{Banach}) first proved that each centrally symmetric convex curve has a conjugate pair of diameters.

\begin{proposition}\label{conjdiam}
The unit circle of any Minkowski plane has a pair of conjugate diameters.
These diameters may be chosen such that their endpoints are extreme points of the unit disc.
\end{proposition}
\proof
By compactness there is a triangle $\triangle\vo\vp\vq$ of maximum area, where $\vp,\vq$ range over all unit vectors.
By this maximum area property it follows that the line through $\vp$ parallel to $\vq$ supports the unit circle, i.e., $\vp\normal\vq$, and vice versa.
If $\vp$ is in the relative interior of the segment $\Segment{\va\vb}$ on the unit circle, then either $\triangle\vo\va\vq$ or $\triangle\vo\vb\vq$ must also have maximum area.
\qed

\begin{figure}
\begin{center}
\includegraphics{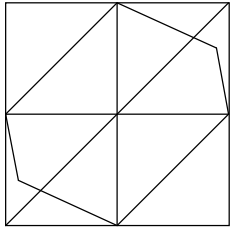}
\end{center}
\caption{Proof of Proposition~\ref{longsegmentprop}}\label{fig23}
\end{figure}

The following proof is a typical application of Proposition~\ref{conjdiam}.
For another, see \cite{MR99k:52028}.
\proof[Proof of Proposition~\ref{longsegmentprop}]
Consider a coordinate system determined by conjugate diameters of the unit circle whose endpoints are extreme points.
Assume that there are at least three pairs of segments of length at least $1$ on the unit circle.
Then there are at least two such segments in some quadrant, say the first.
See Figure~\ref{fig23}.
Let $d$ be the length of the arc of the unit circle lying in the first quadrant.
Then $d\geq 2$ (since each segment is of length at least $1$).
But also $d\leq 2$ by Theorem~\ref{convexcurveineq}.
Thus $d=2$, and it follows that the two segments in the first quadrant share an endpoint, with the other endpoints being the standard unit vectors $\ve_1,\ve_2$, and both are of length exactly $1$, and $\Segment{\ve_2(-\ve_1)}$ is a segment on the unit circle by equality in the Monotonicity Lemma.
It follows that the unit circle is determined.

In particular, the endpoint shared by the two segments in the first quadrant must be in the direction $\frac{1}{2}(\ve_1+\ve_2)$, and there are exactly three pairs of segments of length $\geq 1$.
\qed

Another way of finding conjugate diameters is by considering a parallelogram of minimum area circumscribed to the unit circle.
By minimality of area it is easily seen that the midpoints of the sides of the parallelogram must be on the unit circle, and since the parallelogram is circumscribed, the segments joining midpoints of opposite sides are conjugate diameters.
Funk \cite[p.~92]{Funk} mentions the idea of considering a smallest circumscribed parallelogram, but it is not clear that he knew of the above result.
The earliest reference to this proof using circumscribed parallelograms is Day \cite{MR9:246h}, who also gave a generalization to circumscribed $2n$-gons and to higher dimensions (Section~\ref{higherdim}).
S\"uss \cite{MR16:1046b} gave a proof using the inscribed parallelogram of maximum area, although this is later than Taylor \cite{MR8:588c} and Day \cite{MR9:246h}, who both already cover the higher-dimensional case.

Lenz \cite{MR21:2204} shows that these two ways always give different pairs of conjugate diameters, except if the unit circle is a Radon curve, in which case any diameter has a conjugate diameter.
Auerbach \cite{Auerbach2} also announces that a planar symmetric convex disc always has at least two different pairs of conjugate diameters, but without proof.
Laugwitz \cite{MR16:613d} uses a parametrization of convex curves due to Radon \cite{Radon, MR89i:01142a} to prove that there are at least two conjugate diameters.
Inzinger \cite{MR10:60c} gives a sufficient condition for a symmetric disc to have exactly two pairs of conjugate diameters: no concentric ellipse should intersect it in more than four diametral pairs of points.
For a comprehensive discussion of conjugate diameters of planar convex curves, see Heil and Krautwald \cite{MR40:3431}.
In particular, they show that a pair of conjugate diameters induced by an inscribed parallelogram of largest area corresponds to a pair of conjugate diameters of the dual induced by a circumscribed parallelogram of smallest area, and vice versa.
Weiss \cite{Weiss} proves that by iterating the operation of choosing normals to vectors, one obtains conjugate diameters in the limit.

\subsubsection{Radon curves}\label{radon}
A \emph{Radon curve} is the unit circle of a Minkowski plane in which normality is symmetric, i.e.,
for any $\vx,\vy\neq\vo$, if $\vx\normal\vy$ then $\vy\normal\vx$.
Radon curves were introduced by Radon \cite{Radon, MR89i:01142a}, who gave examples of them and studied some of their properties.
He introduced them concretely, namely as those curves which are affinely equivalent to a curve whose polar is a $90^\circ$ rotation of the original curve.
He showed that any curve for which normality is symmetric must be of this form.
What is essentially needed for the proof is Proposition~\ref{uniquebody}.
He indicated that there exist algebraic curves that are Radon curves, a statement worked out in detail by Leichtweiss \cite{MR33:3191}.
He also showed that these norms are exactly those in which
\[ \length{\va\vb} = \lim_{\epsi\to 0} \frac{\Area(\Segment{\va\vb}+\epsi B)}{2\epsi} .\]
Birkhoff \cite{Birkhoff} independently discovers and constructs Radon curves.
James \cite{MR9:42d} also constructs these curves and gives the example of a two-dimensional norm which is the $p$-norm in the first and third quadrants and the $q$-norm in the second and fourth quadrants, where $1/p+1/q=1$.
Day \cite{MR9:192c} describes the construction in detail, as well as the proof that this construction exhausts all norms for which normality is symmetric (this follows immediately from Proposition~\ref{uniquebody}).
Heil \cite{MR37:5796} notices that the regular $2n$-gons are Radon curves iff $n$ is odd.
Note also that any Radon curve can be approximated by a sequence of polygonal Radon curves, as noted by Day \cite[pp.~332--333]{MR9:192c}.
Textbooks which discuss Radon curves are Busemann \cite{MR17:779a}, Benson \cite{MR35:844}, Leichtweiss \cite{MR81b:52001} and Thompson \cite{MR97f:52001}.
See also Krautwald \cite{Krautwald} for a discussion of Radon curves.

Phelps \cite{MR19:432a} relates symmetry of normality to nonexpansiveness of nearest-point maps.
Let $M$ be a Minkowski space and $S$ a subset.
Then a mapping $\vf:M\to S$ is a \emph{nearest-point map} if $\vf(\vx)$ is a point in $S$ nearest to $\vx$, among all points in $S$.
Phelps shows that a Minkowski plane $M$ is strictly convex and of Radon type iff for any closed convex subset $S$ and any nearest point map $\vf:M\to S$ is \emph{non-expansive}, i.e., for all $\vx,\vy\in M$, $\length{\vf(\vx)\vf(\vy)}\leq\length{\vx\vy}$.
Thus in higher dimensions only Euclidean space has this property.
Related to this are results of Sch\"affer \cite{MR31:3834} and De Figueiredo and Karlovitz \cite{MR35:2130} on the so-called radial projection (see below).

Note that Radon curves may be characterized as those unit circles in which each boundary point is the midpoint of a side of a circumscribed parallelogram of smallest area.
It is necessary to require that each point is a midpoint of a circumscribed parallelogram, as the example of a regular octagon shows.
See also Section~\ref{circumferencesect} where equiframed curves are discussed.

Lenz \cite{MR21:2204} studies various extremal properties of Radon curves.
Let $A_1$ be the maximum area of a quadrilateral inscribed in the unit circle, $A_2$ the minimum area of a circumscribed quadrilateral, and $A_3$ the minimum area of a circumscribed parallelogram.
Lenz shows that $A_3/A_1\leq 2$, with equality only for Radon curves.
It follows from this inequality that whenever the unit circle is not a Radon curve, then the conjugate diameters induced from an inscribed parallelogram of largest area, and from a circumscribed parallelogram of smallest area, must be different (as noted above).
He also shows that $\Area(B)/A_1\leq \Pi(M)/2$, with equality iff $B$ is a Radon curve, that $A_3/\Area(B)\leq 4/3$ if $B$ is a Radon curve (proved earlier for any Minkowski plane by Petty \cite{MR18:760e}), and that $3\leq\Pi(M)\leq\pi$ for Radon curves.
This last inequality was rederived by different methods in \cite{MR37:5796}, and was also rediscovered by Yaglom \cite{MR58:30739}.
See Section~\ref{circumferencesect} for a further discussion of these three results.
As noted in \cite{MR40:3431}, the above inequalities remain true with $A_3$ replaced by $A_2$, since by a result of Dowker \cite{MR5:153m}, if a convex disc is centrally symmetric, then among circumscribed $2n$-gons of smallest area there is a symmetric one (with the same centre as the convex disc).
Lenz also characterizes the Euclidean plane as a Radon plane in which each pair of conjugate diameters partitions the unit circle into four pieces of equal area (or of equal length).

Laugwitz \cite{MR16:613d} characterizes Radon curves as those for which the tangential curvature of the unit circle parametrized by arc length is constant (and then necessarily $0$).

Busemann \cite{MR9:372h}, in solving the isoperimetric problem for Minkowski planes, shows that a plane is of Radon type iff isoperimetrices are circles.

Busemann \cite{MR17:779a} proves that in Radon planes the following definition of area of a triangle is independent of which side is chosen as base: $\Area = hb/2$, where $h$ is the shortest distance from an opposite vertex to a side and $b$ is the length of the side of the triangle.
Tam\'assy \cite{MR23:A4052, Tamassy} shows the converse: if the definition is independent of the choice of side of the triangle, then the plane is of Radon type.

S\"uss \cite{MR12:46d} gives a characterization of Radon curves in terms of chords.
Porcu \cite{MR30:2400} and Heil \cite{MR37:5796} discuss various affine properties of Radon curves.

In \cite{MR31:3834} Sch\"affer introduces a constant, later named the \emph{Sch\"affer constant}, as
\[ \mu := \liminf_{\length{\un{\vx}\un{\vy}}\to 0}\frac{\length{\vx\vy}}{\length{\un{\vx}\un{\vy}}\max(\norm{\vx},\norm{\vy})}.\]
This constant always satisfies $1/2\leq\mu\leq 1$, and Sch\"affer shows that $\mu=1$ iff normality is symmetric, hence in the two-dimensional case, iff the norm is of Radon type.
It therefore follows that the \emph{radial projection} $R:M\to M$, defined by
\[ R\vx := \left\{\begin{array}{lcl} 
              \vx      & \mbox{ if } \norm{\vx}\leq 1 \\ 
              \un{\vx} & \mbox{ if } \norm{\vx}>1,
           \end{array}\right. \]
is non-expansive, a fact also proved by De Figueiredo and Karlovitz \cite{MR35:2130}, and generalized by Karlovitz \cite{MR46:7869} and Gruber \cite{MR51:13648, MR57:1279, MR86k:41032}.
Recently this result was rediscovered in \cite{MR96f:46022}.
Sch\"affer also shows that the \emph{normalization} $\un{\gamma}$ of a curve $\gamma$, defined by $\un{\gamma}(\vx):=\un{\gamma(\vx)}$, has length $\length{\un{\gamma}}\leq\length{\gamma}$.

We now give elementary proofs of propositions in \cite{MR35:2130}, indicating how one can use geometric arguments very reminiscent of Euclidean geometry in proving results on Radon planes.

\begin{proposition}\label{radonprop}
In a Minkowski plane $M$ the following four properties are equivalent.
\begin{enumerate}
\item The plane is a Radon plane.\label{stA}
\item The radial projection is non-expansive.\label{stB}
\item For any two unit vectors $\vx,\vy$, there exists a point $\vp$ on $\Line{\vo\vy}$ nearest to $\vx$ satisfying $\vp\in B$.\label{stC}
\item For any two unit vectors $\vx,\vy$, all points $\vp$ on $\Line{\vo\vy}$ nearest to $\vx$ satisfy $\vp\in B$.\label{stD}
\end{enumerate}
\end{proposition}

\begin{figure}
\begin{center}
\includegraphics{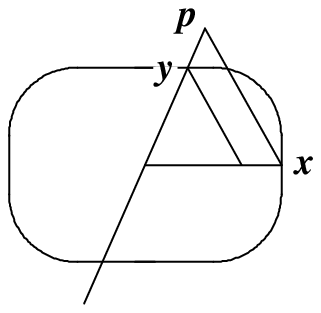}
\end{center}
\caption{\ref{stA}$\Rightarrow$\ref{stD} in Proposition~\ref{radonprop}}\label{fig24}
\end{figure}

\proof
\ref{stA}$\Rightarrow$\ref{stD}.
Let $\vx,\vy$ be unit vectors and $\vp$ a point on $\Line{\vo\vy}$ nearest to $\vx$.
Without loss of generality assume that $\vc\neq\vy$ and that $\vp$ is a positive multiple of $\vy$.
Then $\vp-\vx\normal\vy$, hence $\vy\normal\vp-\vx$.
If $\vp\not\in B$, it now follows that the line through $\vy$ parallel to $\vp-\vx$ intersects $\interior B$, a contradiction.
See Figure~\ref{fig24}.

\ref{stC}$\Rightarrow$\ref{stA}.
Choose any unit vectors $\vx,\vy$ such that $\vx\normal\vy$.

\begin{figure}
\begin{center}
\includegraphics{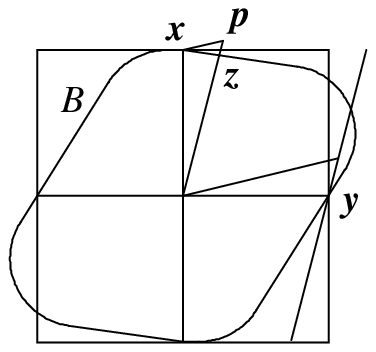}
\end{center}
\caption{\ref{stC}$\Rightarrow$\ref{stA} in Proposition~\ref{radonprop}}\label{fig25}
\end{figure}

Suppose $\vy\notnormal\vx$.
Thus the line through $\vy$ parallel to $\vx$ intersects $\interior B$.
Assume without loss of generality that this is on the same side of $\Line{\vo\vy}$ as $\vx$.
Choose a unit vector $\vz$ on the same side of $\Line{\vo\vx}$ as $\vy$ such that the line $\ell$ through $\vy$ parallel to $\vz$ still intersects $\interior B$.
Then no point on $\ell$ at a shortest distance to $\vo$ can lie on $\Line{\vo\vy}$ or in the open half plane bounded by $\Line{\vo\vy}$ opposite $\vx$.
It follows that all points $\vp$ on $\Line{\vo\vz}$ nearest to $\vx$ are such that $\vz$ is between $\vo$ and $\vp$, a contradiction.
See Figure~\ref{fig25}.

\ref{stB}$\Rightarrow$\ref{stC}.
Let unit vectors $\vx,\vy$ be given.
We are given that for any $\lambda\geq 1$, $\norm{\vx-\lambda\vy}\geq\norm{\vx-\vy}$ and $\norm{\vx+\lambda\vy}\geq\norm{\vx-\vy}$.
By Lemma~\ref{maxineq}, $\lambda\mapsto\norm{\vx-\lambda\vy}$ attains a minimum for $-1\leq\lambda\leq 1$.

\ref{stC}$\Rightarrow$\ref{stB}.
It is sufficient to show that if $\lambda>1$ and $\norm{\vx}\leq 1=\norm{\vy}$, then $\norm{\vx-\lambda\vy}\geq\norm{\vx-\vy}$.
It is given that there exists $\mu$ such that $-1\leq\mu\leq 1$ and $\norm{\vx-\mu\vy}\leq\norm{\vx-\lambda\vy}$.
By Lemma~\ref{maxineq}, $\norm{\vx-\vy}\leq\max(\norm{\vx-\mu\vy},\norm{\vx-\lambda\vy})=\norm{\vx-\lambda\vy}$.
\qed

The implication $\ref{stC}\Rightarrow\ref{stD}$ was already noticed by Amir \cite[\S18]{MR88m:46001}.
Amir implicitly discusses Radon planes in Lemma~10.3 and its corollaries and in \S18 of \cite{MR88m:46001}.

Thele \cite{MR48:6892} studies the Lipschitz constant of the radial projection,
\[\inf\{c>0: \length{R\vx\,R\vy} \leq c\length{\vx\vy} \}. \]
Desbiens \cite{MR94a:46018} shows that this constant and Sch\"affer's constant are in fact equal.
Franchetti \cite{MR82c:46017} shows that the Thele constant of a Minkowski space is equal to the Thele constant of its dual.

Sperner \cite{MR82f:52004} characterizes Radon curves among all  self-dual Minkowski planes in terms of a so-called canonical isomorphism between a self-dual space and its dual introduced by Leichtweiss \cite{MR22:174}.
Alonso and Ben{\'\i}tez \cite{MR91e:46021b} give characterizations of Radon curves in terms of the equivalence of various definitions of orthogonality.

Finally, Radon curves also make an appearance in hyperbolic geometry.
Pinkall \cite{MR85c:52019} characterizes horocyclically convex sets (for any two points in the set the horocyclic segment joining them is also in the set) of constant width in terms of Radon curves.

\subsection{Higher dimensions}\label{higherdim}
We define a set of $d$ diameters of the unit ball of a $d$-dimensional Minkowski space to be \emph{conjugate diameters} if their normalized direction vectors $\vx_1,\dots,\vx_d$ have the property that each $\vx_i$ is normal to each vector in the linear span of the remaining direction vectors.
An \emph{Auerbach basis} of a Minkowski space is such a set of direction vectors.
Note that it is not sufficient to require that $\vx_1,\dots,\vx_d$ are mutually normal, as the following three-dimensional example shows:
Let $\ve_1,\ve_2,\ve_3$ be linearly independent vectors in $\R^3$.
Then $\ve_1,\ve_2,\ve_3$ are mutually normal in the Minkowski space with unit ball $\conv\{\pm\ve_1,\pm\ve_2,\pm\ve_3,\pm r(\ve_1+\ve_2+\ve_3)\}$ for any $r>1$, although they do not form an Auerbach basis.

If $\vx_1,\dots,\vx_d$ is a normalized Auerbach basis, then a \emph{dual normalized basis} is $\phi_1,\dots,\phi_d$, where $\phi_i$ is a unit functional with kernel the hyperplane spanned by $\{\vx_j:j\neq i\}$.
The sequence of pairs $(\vx_i,\phi_i)$ is also called a \emph{biorthogonal system}.
The existence of an Auerbach basis means metrically that for any norm there is a coordinatization of $d$-space such that the norm is majorized by the $1$-norm, and minorized by the $\infty$-norm.

Banach \cite[p.~106 and p.~238]{Banach} defines the notion of a biorthogonal system and says that the existence of an Auerbach basis is due to Auerbach.
However, it seems that Auerbach did not publish anything on the higher-dimensional case \cite{Auerbach1, Auerbach2}.
The first published proofs of existence seem to be Taylor \cite{MR8:588c} and Day \cite{MR9:246h}.

As in the two-dimensional case we have that any Minkowski space has at least two Auerbach bases.
One is induced by a cross-polytope inscribed in the unit ball of maximum volume (Taylor \cite{MR8:588c}, Ruston \cite{MR29:2630}), and the other by the midpoints of the facets of a circumscribed parallelotope of minimum volume (Day \cite{MR9:246h}, Lenz \cite{MR19:977d}).
As noticed by Knowles and Cook \cite{MR50:2879}, these two ways of finding Auerbach bases are dual in the sense that if an Auerbach basis is induced by an inscribed cross-polytope of maximum volume, then any dual basis is induced by a circumscribed parallelotope of minimum volume, and vice versa (the two-dimensional case being noticed in \cite{MR40:3431}; see Section~\ref{conjdiamsection}).

If any minimum volume base and maximum volume base coincide, then by the following result of Lenz \cite{MR19:977d} we have that the space is Euclidean:
Let $V_1$ be the volume of a cross-polytope of maximum volume inscribed in the unit ball, and $V_2$ the volume of a parallelotope of minimum volume circumscribed to the unit ball.
Then $V_2/V_1\leq d!$, with equality iff the space is Euclidean (where the dimension $d\geq 3$).
Lenz proves the characterization under the assumption that the space is smooth or strictly convex, but notices that his argument for the non-uniqueness of Auerbach bases does not need any such assumptions.
Also, it is easily seen that the inequality does not depend on any assumptions of smoothness or strict convexity, as it immediately follows from the fact that a maximal inscribed cross-polytope provides a circumscribed parallelotope of volume $d!$ times that of the cross-polytope.
Plichko \cite{MR96k:46020} removes the assumptions from the above characterization of Euclidean space by showing its equivalence to symmetry of normality.

Plichko \cite{MR96k:46020} also proves that if for any two Auerbach bases of a Minkowski space there is an isometry of the space taking the one Auerbach basis to the other, then the space is Euclidean.

\subsection{Non-symmetric curves}
We briefly remark on the non-symmetric generalization of Radon curves.
Blaschke \cite{MR18:922c} introduced \emph{P-curves}, defined to be the closed convex curves which have a continuous family of inscribed quadrilaterals of maximum area.
He characterizes them as the convex curves for which each diameter has a conjugate.
He shows that Radon curves are exactly the P-curves with a centre of symmetry, and that the difference body of a P-curve is always a Radon curve.
Lenz \cite{MR21:2204} shows that the inequality $A_2/A_1\leq 2$, with $A_2$ the area of a smallest circumscribed parallelogram, and $A_1$ the area of a largest inscribed quadrilateral, is also true for convex curves that are not necessarily centrally symmetric, and that equality holds exactly for P-curves.
He also shows that a P-curve is a curve of constant width with respect to its difference body as unit ball, and that the circumference/width ratio of a P-curve is $\leq \pi$, where the circumference is measured with respect to its central symmetrization, with equality only for those curves which are also affinely equivalent to Euclidean curves of constant width.
He observes that among all P-curves with the same central symmetrization $S$ and width 2, $S$ has the largest area (since central symmetrization of a non-symmetric disc increases area).
Heil and Krautwald \cite{MR40:3431} observe that regular $n$-gons, with $n$ odd, as well as Euclidean curves of constant width, are P-curves.
Porcu \cite{MR30:2400} also proves some properties of P-curves.
Martini \cite{MR91k:52018} gives a characterization of P-curves (see the Concluding Remarks at the end of \cite{MR91k:52018}).

In higher dimensions, Krautwald \cite{MR82j:52012} gives non-symmetric analogues of the characterization of ellipsoids due to Lenz \cite{MR19:977d} mentioned above in Section~\ref{higherdim}: In a space of dimension $d\geq 3$, a convex body is an affine image of a Euclidean body of constant width iff the vertices of an inscribed cross-polytope of maximum volume are on the faces of a circumscribed parallelotope of minimum volume, iff the ratio between the volumes of these two polytopes is $d!$.

\section{Minkowski circles}
\subsection{Circumscribed and inscribed circles}\label{circumsect}
As mentioned earlier, a Minkowski plane is strictly convex iff there is \emph{at most} one circumscribed circle through any three non-collinear points.
We also have the following

\begin{proposition}
A Minkowski plane is smooth iff through any three non-collinear points there is \emph{at least} one circumscribed circle.
\end{proposition}

Sufficiency of smoothness in the plane was proved by Mayer \cite{Mayer}, and even earlier by Zindler \cite[Satz 39, p.~53]{ZindlerII} for the special case of equilateral triangles (which is already general if one uses affine transformations).
We here sketch the proof of Kramer and N\'emeth \cite{MR57:7384, MR50:14504} which, although based on Brouwer's fixed point theorem, is very clear and is immediately generalizable to higher dimensions.

\proof
Let $\vp_i$, $i=1,2,3$, be the three non-collinear points.
Let $\vc$ be their centroid, and let $\vp_i'=\un{\vp_i-\vc}$ be the directions from the centroid to the points.
For any point $\vx$ in the unit ball $B$, let $\vq_i(\vx)$ be the point of $B$ on the ray with origin $\vx$ and direction $\vp_i'$ furthest from $\vx$.

\begin{figure}
\begin{center}
\includegraphics{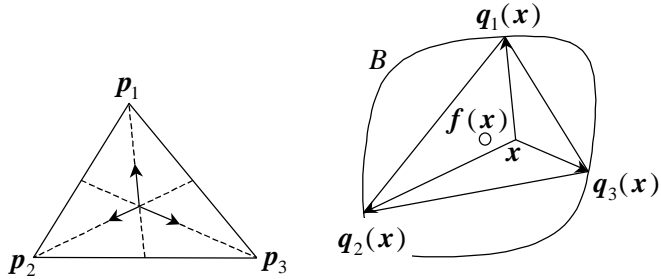}
\end{center}
\caption{Defining $\vf$}\label{fig26}
\end{figure}

\begin{figure}[b]
\begin{center}
\includegraphics{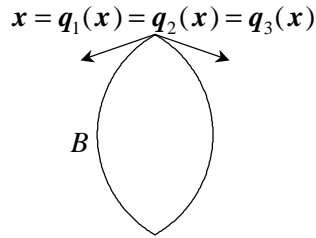}
\end{center}
\caption{Since $B$ is smooth, the fixed point is not on the boundary}\label{fig27}
\end{figure}

It is easily seen that $\vq_i$ are continuous mappings (even if $B$ is not strictly convex).
Define $\vf:B\to B$ by $\vf=\frac{1}{3}(\vq_1+\vq_2+\vq_3)$.
See Figure~\ref{fig26}.
By the Brouwer fixed point theorem, $\vf$ has a fixed point $\vx_0$.
If $\vx_0$ is on the boundary of $B$, then it is easily seen that $B$ is not smooth (see Figure~\ref{fig27}).
Note that it is only at this point that we need smoothness.

Thus $\vx_0\in\interior B$, and it is easily seen that the points $\vp_i(\vx_0)$ are positively homothetic to $\vp_i$, which is equivalent to the fact that $\vx_i$ have a circumscribed circle.

\begin{figure}
\begin{center}
\includegraphics{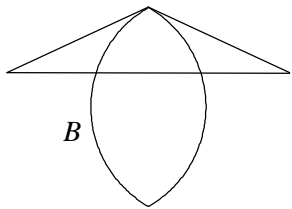}
\end{center}
\caption{A triangle without a circumscribed circle}\label{fig28}
\end{figure}

The converse is clear from Figure~\ref{fig28}.
\qed

The above proof even generalizes to the non-centrally symmetric case, which was, according to Kramer and N\'emeth, a conjecture of Tur\'an.
Makeev also proved this result somewhat later \cite{MR88m:55001}, using more complicated topological methods.
See also \cite{MR95d:52006}.

By contrast any triangle in a Minkowski plane has a unique inscribed circle.
This can be proved using Glogovskii's definition of angular bisectors \cite{MR45:2597} exactly as in the Euclidean plane:
For any angle with sides $\Ray{\va\vb}$ and $\Ray{\va\vc}$, there exists a unique ray $\Ray{\va\vd}$ with the property that each point on $\Ray{\va\vd}$ is equidistant to $\Ray{\va\vb}$ and $\Ray{\va\vc}$.
(This is true even if the space is not smooth.)
If we call $\Ray{\va\vd}$ the \emph{angular bisector} of angle $\myangle\vb\va\vc$, then the three angular bisectors of a triangle intersects in the centre of the unique inscribed circle of the triangle.
This observation in the case of the Minkowski plane with a regular hexagon as unit circle, was made by Sowell \cite{MR91a:51003}.

\subsection{Equilateral $n$-gons inscribed in the unit circle}
In \cite{MR93e:52038} it is shown that for any $n\geq 3$ there is an equilateral $n$-gon inscribed in the unit circle, starting at any given $\vx$.
For a fixed $\vx$, although the $n$-gon need not be unique, its side length is unique.

This result leads to four extremal problems: for given $n$, find an inscribed equilateral $n$-gon of smallest or largest side length, and minimize or maximize this quantity over all Minkowski planes.
In \cite{MR93e:52038} the largest side length is considered, and it is shown that for $n\leq 6$ (when this side length is at least $1$) finding this side length is equivalent to finding the smallest circle containing a packing of $n$ unit discs, essentially by using Lemma~\ref{packingineq}.
They find that if $\delta_M(n)$ is the largest side length in the Minkowski plane $M$, then $\delta_M(4)\geq \sqrt{2}$, with equality for the ellipse.
There are also other planes where equality holds, see Gao \cite{MR98a:46031}.
In \cite{MR93e:52038} it is also conjectured that $\delta_M(3)\geq 1+1/\sqrt{2}$, equality holding if the unit ball is a regular octagon.
Linhart \cite{MR83d:52015} makes an equivalent conjecture.
Doliwka and Lassak \cite{MR96d:52002} show that $\delta_M(5)\geq\sqrt{5}-1$, with equality if e.g.\ the unit circle is an affine regular decagon.

\subsection{Characterization of the Minkowski circle}
Valentine \cite{MR13:377g} gives the following characterization of discs in the Minkowski plane.

\begin{proposition}
A bounded closed subset $S$ of a Minkowski plane is a Minkowski disc iff for any distinct $\vx,\vy\in S$ there exists a semicircle in $S$ with diameter $\Segment{\vx\vy}$.
\end{proposition}

Valentine also shows that the above is true for Minkowski spaces if ``semicircle'' is replaced by ``hemisphere''.

W.-T.\ Hsiang \cite{Hsiang} shows (for the case of the Euclidean plane) the following

\begin{proposition}
A Jordan curve $\gamma$ in a Minkowski plane is a Minkowski circle iff for any distinct $\vx,\vy$, both in the interior or both in the exterior of $\gamma$, there exists a circle through $\vx$ and $\vy$ not intersecting $\gamma$.
\end{proposition}

This is a special case of the following result of Petty and Crotty \cite{MR41:2538}.

\begin{proposition}
A non-empty bounded open subset $S$ of a Minkowski space is an open ball iff for any distinct $\vx,\vy\in S$ there exists an open ball in $S$ with $\vx$ and $\vy$ as boundary points.
\end{proposition}

The above result is in fact proved for a more general class of metric spaces in \cite{MR41:2538}.
See also \cite{Israel}, where essentially the same result is posed as a problem (for Banach spaces) and solved.

See O. Haupt \cite{Haupt} and also \cite{Mayer} for a study of certain sets satisfying a property related to the condition of Valentine.

Petty and Crotty also give the following characterization of balls in Minkowski spaces: a bounded open convex subset of a Minkowski space is an open ball iff it has constant width and possesses an equichordal point (i.e.\ all chords through this point are of equal length).

The following characterization was proved by Groemer \cite{MR25:2515} for the Euclidean case.

\begin{proposition}
Let $C$ be a convex body in a Minkowski space with the property that there exists a point $\vp$ such that for any boundary point $\vq$ of $C$ there is a supporting hyperplane $H$ at $\vq$ such that the unit ball has a supporting hyperplane at $\un{\vq-\vp}$ parallel to $H$.
Then $C$ is a ball.
\end{proposition}

The above proposition follows immediately from the following characterization of homothets of convex bodies (proved by Sch\"affer \cite[4A]{MR57:7120}).

\begin{proposition}\label{uniquebody}
Let $K_1$ and $K_2$ be convex bodies in $d$-dimensional space both containing $\vo$ as interior point.
Assume that for each ray with origin $\vo$ there are parallel hyperplanes supporting $K_1$ and $K_2$ at the points where the ray meets the boundaries of $K_1$ and $K_2$.
Then $K_1=\lambda K_2$ for some $\lambda>0$.
\end{proposition}

Witsenhausen \cite{MR45:5939} proved the following

\begin{proposition}
Let $\gamma$ be a simple closed rectifiable curve in a Minkowski plane $M$.
For each $\vx$ on the curve let $\vx'$ be the point whose distance along $\gamma$ is half its length $\length{\gamma}$.
Then \[\length{\gamma}\geq \Pi(M) \min_{\vx\in\gamma} \length{\vx\vx'}. \]
\end{proposition}
Note that Witsenhausen proved a similar result in any Minkowski space, where $\Pi$ is replaced by half the \emph{girth} of the unit ball (see Sch\"affer \cite{MR36:1959}).
It is known that in the Euclidean plane equality in the above inequality characterizes the circle \cite{MR48:12290}, but Witsenhausen notes that for example in the rectilinear plane there are other curves than circles for which equality holds.
However, it is probable that in strictly convex Minkowski planes equality characterizes circles.

The results of Goodey and Woodcock (see Section~\ref{intercircsect}) may also be considered as a characterization of the Minkowski circle.

\subsection{Circumference and area of the unit circle}\label{circumferencesect}
As will be seen below, it is profitable to study the circumference of the unit circle together with the area of the unit circle (where there are various intrinsic ways of fixing the unit area).
In the first case we use area to study circumference, and in the second case we use circumference to study area.

Go{\l}ab \cite{Golab} proved that the circumference of the unit circle is at
least $6$, with equality iff the unit circle is an affine regular
hexagon, and at most $8$, with equality iff the unit circle is a
parallelogram.
This was rediscovered by later authors, e.g., Petty \cite{MR18:760e}, Sch\"affer \cite{MR36:1959}, Laugwitz \cite{MR16:613d}, Re\v setnyak \cite{MR15:819d}.
Sch{\"a}ffer \cite{MR47:5732} shows that the circumference is equal to that of the unit circle of the dual plane; see also Thompson \cite{MR52:4138}.
This is posed as a question for the $p$-norms as recently as \cite{AdlerTanton}.
See \cite{MR97f:52001} for an exposition of these results.

There are many papers on the circumference of unit circles of non-symmetric norms; see \cite{MR30:2396, MR34:1925, MR33:3188, MR40:3434, MR49:3687, MR81e:52002, MR99c:52009}.
See also Heil \cite{MR37:5796} for inequalities relating cirumference and area, deduced using inequalities from convex geometry.
Recent references on the circumference in the symmetric case are the following.

Franchetti and Votruba \cite{MR57:10594} show that four times the absolutely summing constant of a Minkowski plane equals the circumference of its unit circle.
Chalmers, Franchetti and Giaquinta \cite{MR96m:46018} show that in a \emph{symmetric} Minkowski plane (for some basis $\vx,\vy$ the norm satisfies $\norm{\abs{\lambda}\vx+\abs{\mu}\vy}=\norm{\lambda\vx+\mu\vy}=\norm{\mu\vx+\lambda\vy}$), the circumference of the unit circle is at least $2\pi$, with equality iff the plane is Euclidean.
See Adler and Tanton \cite{AdlerTanton} and Euler and Sadek \cite{EulerSadek} for estimates of the numerical values in the $\ell_p$ planes.
Ghandehari and Pfiefer \cite{Ghandehari} give formulas for the circumference of the unit circle if it is a regular $2n$-gon.

We now state a result that is not so well known, due to
Lenz \cite{MR21:2204} and rediscovered by Yaglom \cite{MR47:4147}.
See also Heil where a different proof is given.
We follow the proof given by Lenz.
\begin{proposition}\label{Radoncircumprop}
In a Radon plane, the circumference of the unit circle is at least
$6$ with equality iff the unit circle is an affine regular hexagon, and at most $2\pi$, with equality iff the plane is Euclidean.
\end{proposition}

Its proof is given after we discuss a notion of area for Minkowski planes.

There are many ways of defining unit area. One way is to consider the parallelogram of maximum area inscribed in the unit circle, and let the area of this parallelogram be $2$.
This gives the usual area in the Euclidean plane.
Then the area of the unit circle is at least $2$ (equality iff the plane is rectilinear --- trivial), and at most $\pi$ (equality iff the plane is Euclidean --- see Sas \cite{Sas} and Macbeath \cite{MR12:526e} for a simpler proof).
Instead of inscribed parallelograms, we may take an inscribed $n$-gon of maximum area, and normalize for the Euclidean plane.
Note that, since the unit circle is centrally symmetric, by a theorem of Dowker \cite{MR5:153m} there exists a centrally symmetric $n$-gon achieving the maximum if $n$ is even.
Again, we will have the Euclidean plane as the upper extreme by the result of Sas.
The lower extreme is again trivial if $n$ is even (attained iff the unit ball is a centrally symmetric $n$-gon), but is still interesting for $n$ odd.

With this definition of area, the area/circumference ratio is \emph{a priori} in the interval $[1/4,\pi/6]$, with the lower bound best possible.
The upper bound can be improved to $1/2$, with equality iff the unit circle is a Radon curve (Lenz \cite{MR21:2204} and Yaglom \cite{MR47:4147}).

\begin{proposition}\label{areacircumineq1}
In any Minkowski plane, let $P$ be a parallelogram of maximum area
inscribed in the unit circle.
Then $\Area(S)/\Area(P)\leq\length{S}/4$, with equality iff $S$ is a
Radon curve.
\end{proposition}

\proof
By an affine transformation we may assume without loss of generality that the diagonals of $P$ are $\Segment{(-\ve_1)\ve_1}$ and $\Segment{(-\ve_2)\ve_2}$, with $\ve_1$ and $\ve_2$ the standard unit vector basis of $\R^2$.
By Proposition~\ref{conjdiam} the diagonals of this parallelogram are conjugate diameters.
Note that if we take any two unit vectors $\vu,\vv$, then $\abs{\det[\vu,\vv]}\leq 1$, with equality implying that $\vu$ and $\vv$ are on conjugate diameters.

Parametrize $S$ by arclength $\gamma:[0,2\Pi(M)]\to M$.
Let $\vp:[0,2\Pi(M)]\to M$ be the unit (normalized in $M$) right tangent vector.
Then $\vp$ is right continuous, and also left continuous except at countably many points.
We may now find the area of $S$ by integrating:
\begin{eqnarray*}
\Area(S) & = & \int_0^{2\Pi(M)}\frac{1}{2}\det[\gamma(t),\vp(t)] dt \\
& \leq & \frac{1}{2}\length{S}.
\end{eqnarray*}
Equality holds iff $\det[\gamma(t),\vp(t)]=1$ except possibly at nonsmooth points of $S$.
Approximating nonsmooth points from the right ($\vp$ is right continuous), we obtain $\det[\gamma(t),\vp(t)]=1$ for all $t$, hence each diameter is conjugate to some other diameter, and it follows that $S$ is a Radon curve.
\qed

\proof[Proof of Proposition~\ref{Radoncircumprop}]
This follows from the previous proposition, as well as the result of Sas \cite{Sas} that the area of the unit circle is at most $\pi$, with equality iff it is a Euclidean circle.
\qed

In higher dimensions, letting the volume of the inscribed cross-polytope of maximum volume be $2^d/d!$, the volume is at least $2^d/d!$ (again trivial), and at most that of the Euclidean unit sphere (following from the higher dimensional counterpart of the theorem of Sas, due to Macbeath \cite{MR12:526e}).

Another way of defining area is to let the area of the circumscribed parallelogram of least area have area $4$.
Then the area of the unit disc is at most $4$ with equality iff it is a parallelogram (trivial), and at least $3$ with equality iff it is an affine regular hexagon (proved by Lenz \cite{MR21:2204} in the case of Radon curves, and in general by Petty \cite{MR18:760e}, and later also by Babenko \cite{MR89m:52015} and Pe{\l}czynski and Szarek \cite{MR92b:52014}).
We now give a simple proof of this result, by deducing it from a seemingly new area/circumference inequality dual to Proposition~\ref{areacircumineq1} (which Lenz \cite{MR21:2204} missed, although he had all the tools).
Following Pe{\l}czynski and Szarek \cite{MR92b:52014}, we say that a convex body is \emph{equiframed} if each boundary point is contained in the boundary of a circumscribed parallelotope of minimum volume.
Note that an equiframed centrally symmetric convex disc is not necessarily a Radon curve, as the example of a regular octagon shows.
More generally, any regular $2n$-gon is equiframed.
However, it can be shown that strictly convex or smooth equiframed curves must be of Radon type.
 
\begin{proposition}\label{areacircumineq2}
In any Minkowski plane, let $P$ be a parallelogram of minimum area
circumscribed to the unit circle.
Then $\Area(S)/\Area(P)\geq\length{S}/8$, with equality iff $S$ is an equiframed convex disc.
In particular, $\Area(S)/\Area(P)\geq 3/4$, with equality iff $S$ is an affine regular hexagon.
\end{proposition}

\proof
After an affine transformation we may assume that $P=\conv\{\pm\ve_1,\pm\ve_2\}$.
As in the proof of Proposition~\ref{areacircumineq1}, we parametrize $S$ by arc length and use the normalized right tangent vector $\vp$.
Note that $\abs{\det[\gamma(t),\vp(t)]}\geq 1$, with equality implying that $\gamma(t)$ lies on a circumscribed parallelogram of minimum area, namely the parallelogram with a side parallel to $\vp(t)$ and a side parallel to a supporting line of the unit circle at $\un{\vp(t)}$.
We now evaluate the area:
\begin{eqnarray*}
\Area(S) & = & \int_0^{2\Pi(M)}\frac{1}{2}\det[\gamma(t),\vp(t)] dt \\
& \geq & \frac{1}{2}\length{S}.
\end{eqnarray*}
As in the proof of Proposition~\ref{areacircumineq1}, we have equality iff $\det[\gamma(t),\vp(t)]\geq 1$ for all $t$, iff each point is on a circumscribed parallelogram of minimum area.
\qed

Babenko \cite[Proposition 2.2]{MR89m:52015} shows that any unit disc contains a Radon curve for which a minimum circumscribed parallelogram has the same area as a minimum circumscribed parallelogram of the original disc.

See Babenko \cite{MR89m:52015}, Ball \cite{MR90i:52019} and Pe{\l}czynski and Szarek \cite{MR92b:52014} for the higher dimensional question, where there are only partial results.
Ball gives an asymptotic relationship between this ratio and the so-called volume ratio (the fourth way of defining the volume of a unit ball, defined below using the Loewner ellipsoid).

A third way of defining area is to consider the inscribed affine regular hexagon of maximum or of minimum area, and to normalize accordingly.
This is equivalent to asking for the extremes of the ratio of the areas of equilateral triangles to the area of the unit disc.
Reimann \cite{MR89d:52027} proved the inequality below, and the equality cases were characterized by Wellmann and Wernicke \cite{MR93d:51044}.
Wernicke \cite{MR96m:51006} extends these results to Reuleaux triangles.

\begin{proposition}\label{eqtriarea}
Let $\Delta$ be a unit equilateral triangle in the Minkowski plane $M$.
Then
\[ \frac{1}{8} \leq \frac{\Area(\Delta)}{\Area(B)} \leq \frac{1}{6}. \]
There is equality on the left iff $M$ is the rectilinear plane or has a centrally symmetric hexagon (not necessarily affine regular) as unit disc, and equality on the right iff the unit disc is an affine regular hexagon and some translate of $\Delta$ has $\vo$ and two vertices of the unit circle as its own vertices.
\end{proposition}

Weissbach and Wernicke \cite{MR97m:52021} considers this ratio in higher dimensions, and obtain an upper bound of $\binom{2d}{d}^{-1}$, which follows from the Rogers-Shephard inequality \cite{MR19:1073f}.
They show that there is no positive lower bound for $d\geq 3$; this is immediate since in the $d$-dimensional $\ell_\infty$ space one can find $d+1$ equilateral points in a hyperplane, and it is then possible to shift one point, without losing the equilateral property, out of the hyperplane by an arbitrarily small distance, to obtain a simplex of arbitrarily small volume.

A fourth way of defining area is to let the Loewner ellipsoid have area $\pi$.
Then the area of the unit circle is at most $4$, equality if the plane is rectilinear (due to Keith Ball \cite{MR90i:52019}), and at least $\pi$ (trivial).
Again the area/circumference ratio may be examined.

In higher dimensions the volume is at most $2^d$, equality iff the unit ball is a parallelotope (Ball \cite{MR90i:52019}), and at most that of the Euclidean unit sphere.

Finally, see \cite[Chapters 5--7]{MR97f:52001} for various ways of defining area, i.e., $(n-1)$-dimensional content, the two most important definitions coming from Buseman \cite{MR11:200j, MR12:527b} and Holmes and Thompson \cite{MR81k:52023}.
There are various upper and lower bounds for the surface area of the unit ball in dimensions at least three, most of them not sharp -- see \cite{MR81k:52023} and \'Alvarez \cite{Alvarez2}.

\addcontentsline{toc}{section}{References}

\providecommand{\bysame}{\leavevmode\hbox to3em{\hrulefill}\thinspace}

\end{document}